\documentclass{article}
\usepackage{graphicx}
\usepackage{pstricks, pst-node}
\usepackage{natbib}
\usepackage[font=small,labelfont=bf]{caption}
\usepackage{subfig}
\usepackage{amsmath}
\usepackage{amsthm, amssymb, amsfonts}
\allowdisplaybreaks[1]
\usepackage{appendix}
\usepackage{booktabs}
\usepackage{ctable, rotating}
\usepackage{multirow}
\usepackage{appendix}
\usepackage[notablist, tablesfirst, nomarkers, figlist ]{endfloat}
\usepackage[finalnew]{trackchanges}
\addeditor{JK}

\usepackage[top=25mm, bottom=30mm, left=30mm, right=20mm]{geometry}
\usepackage{setspace}
\newcommand{\beq}{\begin{equation}}
\newcommand{\eeq}{\end{equation}}
\newcommand{\beqn}{\begin{eqnarray}}
\newcommand{\eeqn}{\end{eqnarray}}

\newcommand{\ql}{\textquoteleft}
\newcommand{\qr}{\textquoteright}

\begin{document}

\title{Point process models for fine-resolution rainfall}
\author{Jo Kaczmarska, Valerie Isham,\\
Department of Statistical Science \\
University College London, Gower Street, London WC1E 6BT, UK\\\\
Christian Onof, \\
Department of Civil and Environmental Engineering \\
Imperial College London, SW7 2AZ, London, UK. }

\date{\today}
\maketitle

\begin{abstract}
In a recent development in the literature, a new temporal rainfall model, based on the Bartlett-Lewis clustering mechanism and intended for sub-hourly application, was introduced.  That model replaced the rectangular rain cells of the original model with finite Poisson processes of instantaneous pulses,  allowing  greater variability in rainfall intensity over short  intervals.   In the present paper, the basic instantaneous pulse model is first extended to allow for randomly varying storm types.  A systematic comparison of a number of key model variants, fitted to 5-minute rainfall data from Germany, then generates further new insights into the models, leading to the development of an additional model extension, which introduces dependence between rainfall intensity and duration in a simple way.   The new model retains the original rectangular cells, previously assumed inappropriate for fine-scale data, obviating the need for the computationally more intensive instantaneous pulse model.

\end{abstract}

\begin{center}
Accepted for publication in {\it Hydrological Sciences Journal}, September 2013.
\end{center}
\newpage

\section{Introduction}
For both operational and design purposes, hydrologists require long series of rainfall. The time-steps required vary from the daily (for large catchment studies) down to a few minutes for small, typically urban catchments. Observed data series are however often too short, particularly at fine time-scales. While there is an abundance of long daily records in the UK, the number of hourly rainfall records longer than 10 years is of the order of 100. Records of 5 or 10 min rainfalls are generally short as they are recorded over the duration of a particular study (e.g. Hyrex experiment, \cite{Moore00}), and thus include only a limited range of rainfall variability.

The possibility of obtaining ensembles of long series of realistic rainfall data at a range of scales has been the motivation for the development of stochastic rainfall generators over the past four decades. The realism of the data is typically measured by the model's ability to reproduce standard statistics of the time-series of rainfall depth (mean, variance, skewness, autocorrelations), the proportion of dry periods, and the extreme behaviour of rainfall depths at different scales of interest. Unsurprisingly, given the complexity and diversity of the precipitation generating mechanisms, the accurate reproduction of so many features has however so far eluded existing rainfall models. Clear progress is however detectable, and this paper is a contribution to one well-established approach to rainfall modelling.

This approach involves the representation of the physical rainfall process in a realistic, if simplified way, such that the hierarchical structure of rainfall is explicitly incorporated, and the parameters have interpretable meanings.  Introduced in two seminal papers  by \cite{Rod87, Rod88}, it applies Poisson cluster processes to the underlying unobserved continuous-time rainfall process with the rainfall contributed by a sequence  of cells within storms; in these models, discrete-time properties are obtained by aggregation and used to fit to discrete observations. \change[JK]{With this approach, it is possible to reproduce properties of rainfall simultaneously at several time-scales, and this is one of its principal advantages.}{One of the principal advantages of this approach is that it makes it possible to reproduce properties of rainfall simultaneously at several time-scales.}

The clustered point process based models have been extended and validated with a range of different types of rainfall (e.g. \cite{Khaliq96, Ver97, Cam00, Smithers02, Vand11}; see reviews \cite{Onof00} and \cite{Whe05}). These studies show the flexibility of this modelling approach as a tool for reproducing standard rainfall statistics at a range of scales from hourly to daily.  However, because these models all assume that the rainfall cells contribute to the total precipitation through a constant rainfall intensity over the life of the cell (rectangular pulse), the models were not considered appropriate for sub-hourly rainfall, which is also a significant requirement, for example for the design of stormwater sewerage systems.  In order to be able simultaneously to represent subhourly rainfall, as well as rainfall at hourly and daily timescales, it was thought necessary to introduce a third level of (sub-cell) temporal structure, and \change[JK]{recent work \mbox{\citep{Cow07}} introduced}{instantaneous pulses of rain were introduced} to address this apparent shortcoming\add[JK]{\mbox{\citep{Cow07, Cow11}}}.  \change[JK]{However, until now no study has investigated the issue of the performance of rectangular-pulse models at sub-hourly time-scales.}{Recently, however, a spatial-temporal rectangular pulse model applied to sub-hourly data \mbox{\citep{Cow10b}} showed a satisfactory performance}.

\change[JK]{The present paper aims}{In the present paper we present a systematic study in order }to clarify the drivers behind model performance at a sub-hourly timescale, and to identify the optimal choice for fine-scale data.  In Section \ref{models}, the models to be compared are summarised.  These include a new instantaneous pulse version which allows for variation between storms in a parsimonious way, following the approach of \cite{Rod88}.  Section \ref{fitting} briefly outlines the fitting methodology.  In Section \ref{comp}, we present a structured comparison of the performance of the basic rectangular pulse model against the instantaneous pulse version, including also two ways of allowing between-storm variation.  In the literature, two types of clustering have been considered.  We focus here solely on the Bartlett-Lewis suite of models, rather than on models based on the Neyman-Scott clustering mechanism, principally because  methodology for an additional level of clustering has been developed for the former, but the two mechanisms generally exhibit similar performance \citep{Rod87}.  Parameter uncertainty is then considered in Section \ref{param.id}, with potential further improvements to the optimal model discussed in Section \ref{improve}, followed by conclusions in Section \ref{conclusions}.

Note that although we are focusing here purely on temporal models i.e. those fitted to a single site, such models may readily be extended to the spatial dimension and fitted to rain-gauges following the approach of \cite{Cow95} or Chapter 5 of \cite{Whe00}.  Also, although the models themselves are stationary, they may be used to produce simulations that allow for climate change as part of a downscaling methodology.  Examples in the existing literature of this type of approach include \cite{Kil07} and \cite{Bur10b}.  A paper showing how the fitting of these models may be adapted to allow for climate change has been submitted.

\section{Specification of the Bartlett-Lewis suite of models}\label{models}

\subsection{Summary of Existing Models}\label{exist}

\begin{figure}
\subfloat[BLRP model]{\includegraphics[trim=4.8cm 15.8cm 5.5cm 4.5cm, clip=true,  width =7.6cm]{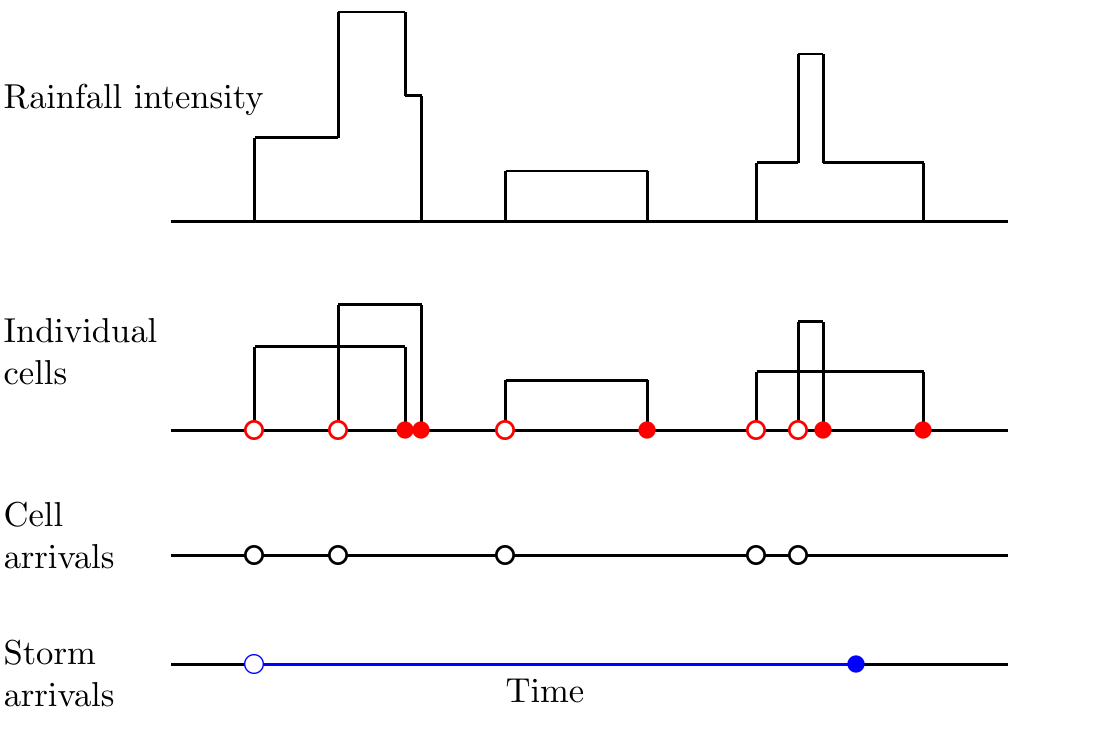}\label{stormpic}}
\hspace{2ex}
\subfloat[BLIP model]{\includegraphics[trim=4.8cm 15.8cm 5.5cm 4.5cm, clip=true,  width =7.6cm]{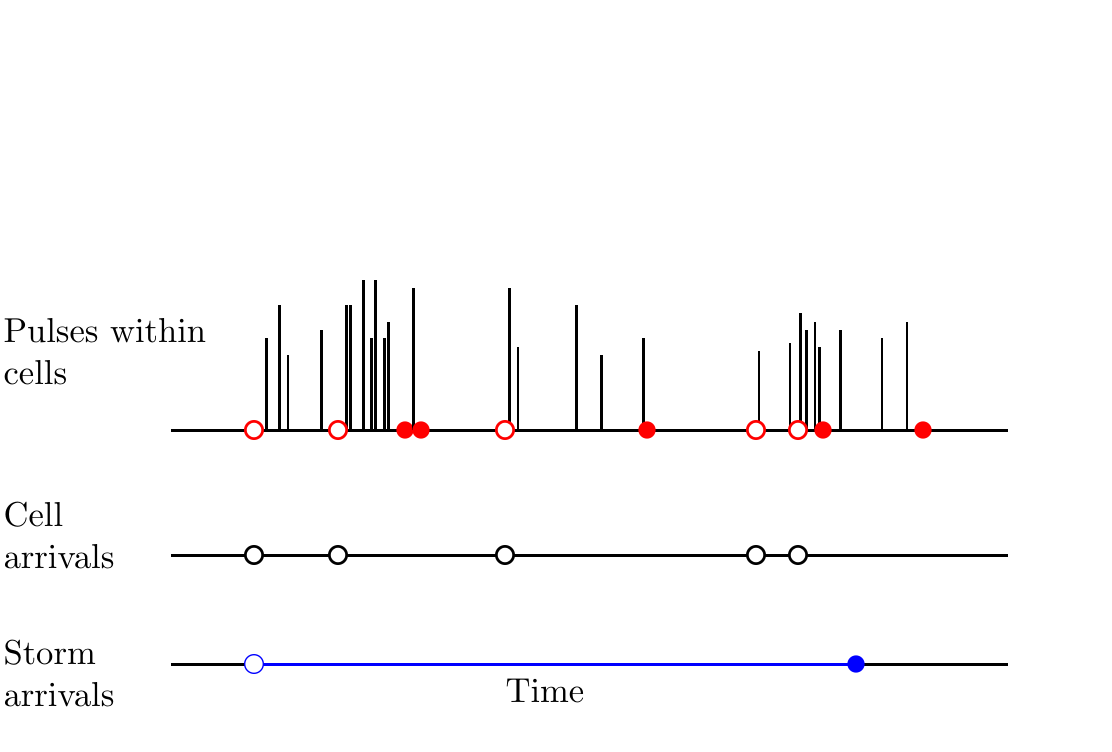}\label{stormpic2}}
\caption{\small{Illustration of a single storm; storm and cell origins are denoted by open circles, and terminations by filled circles. In the BLRP model, each rain cell is assumed to have a constant intensity, whereas in the BLIP model, each cell consists of a series of instantaneous pulses.}}
\end{figure}

In the basic Bartlett-Lewis Rectangular Pulse (BLRP) model, storms arrive in a Poisson process of rate $\lambda$, each storm generating a cluster of cell arrivals. The Bartlett-Lewis clustering mechanism assumes that the time intervals between successive cells are independent, identically distributed random variables  (whereas in the Neyman-Scott model, it is the temporal distances of the cells from their storm origin which are independent and identically distributed).  It is normally assumed that the intervals between cells are exponentially distributed, so that the cell arrivals constitute a secondary Poisson process of rate $\beta$.  Each cell is associated with a rectangular pulse  of rain, of random duration, $L$, and with random intensity, $X$.   In the simplest version of the model, these are both assumed to be exponentially distributed with parameters $\eta$ and $1/\mu_X$ respectively, and are independent of each other. The cell origin process terminates after a time that is also exponentially distributed with rate $\gamma$.  This basic version thus has five parameters in total. Both storms and cells may overlap, and the total intensity of rain at any point in time, $Y(t)$ is given by the sum of all pulses \ql active\qr\, at time $t$.   \add[JK]{ The process in respect of a single storm is illustrated in Figure {\ref{stormpic}}}. Additional flexibility can be added by allowing for a distribution with more parameters for pulse intensities.  A distribution with a longer tail may help in particular with the fit of extreme values, and popular variants include the Gamma, Weibull and Pareto distributions.  One additional parameter is required in order to use either of these.

In order to allow for different types of rainfall, multiple superposed processes can be used \citep{Cow04, Cow07}.  Due to parameter identification issues, the number of processes is typically limited to just two, which can be thought of as representing heavy, short-duration convective and lighter, long-duration stratiform types of rainfall.  However, a \change[JK]{much}{potentially} more parsimonious approach to enable variation between storms is to randomise the cell duration parameter and related temporal storm characteristics.  The Random Parameter Bartlett-Lewis (BLRPR) model \citep{Rod88} extends the basic model by allowing the parameter $\eta$, that specifies the duration of cells, to vary randomly between storms.  This is achieved by assuming that the $\eta$ values for distinct storms are independent, identically distributed random variables from a gamma distribution with index $\alpha$ and rate parameter $\nu$.  The model is re-parameterised so that, rather than keeping the cell arrival rate, $\beta$, and the storm termination rate, $\gamma$ constant for each storm, it is the ratio of both of these parameters to $\eta$ that is kept constant.  Thus, for a higher $\eta$ (i.e. typically shorter cell durations), we have correspondingly shorter storm durations, and shorter cell interarrival times.  Essentially the effect is that all storms have a common structure, but distinct storms occur on different (random) timescales.

The durations of cells and storms (more precisely cell origin processes) are both exponentially distributed, conditional on the cell duration parameter, $\eta$.  Their unconditional distributions are Pareto type II.  This heavy-tailed distribution has an infinite mean if $\alpha$ is less than 1, and an infinite variance if  $\alpha$ is less than 2. Also, in terms of the aggregated rainfall process, it turns out that, for values of $\alpha$ smaller than 3, the variance is infinite, and for values smaller than 4, the skewness is infinite.  This is potentially problematic.  For example, in practice it has been found that simulations with unconstrained values of $\alpha$ occasionally generate unrealistically long periods of rainfall \citep{Onof09, Ver10}.  This can be addressed by setting constraints on $\alpha$, or rejecting storms or cells beyond a certain length or cells with an excessive intensity within any simulations.  Alternatively, the gamma distribution for the cell duration parameter, $\eta$, may be truncated, with support ($\varepsilon, \infty$) \citep{Onof09}. The lower limit, $\varepsilon$, for the integrals over $\eta$ can be pre-specified, or alternatively, can constitute a further parameter to be determined.

The Bartlett Lewis Instantaneous Pulse (BLIP) model \citep{Cow07}, intended for fitting to fine-scale (of the order of five to fifteen minute) data, has a minimum of six parameters (one more than the original Bartlett-Lewis model)\add[JK]{, and is illustrated in Figure {\ref{stormpic2}}}.  As in the BLRP model, storm origins arrive in a Poisson process of rate $\lambda$, and each storm origin initiates a Poisson process of cell origins of rate $\beta$.  In contrast to the basic Bartlett-Lewis  model, however, it is not assumed that there is a cell at the storm origin itself, so a storm may have no rainfall.  This is purely for mathematical convenience and does not lead to any loss of generality.  Each cell origin initiates a further Poisson process of rainfall pulses of rate $\xi$.  Again, it is not assumed that there is a pulse at the cell origin, so a cell may have no rainfall.  Note that the pulses are instantaneous - they have a depth, but no duration.  This Poisson process of instantaneous pulses replaces the rectangular pulse assumption of the original Bartlett-Lewis model.  Both the storm duration (the duration of the cell origin process), and the cell duration are assumed to be exponentially distributed, the former with rate $\gamma$, and the latter with rate $\eta$.  The process of pulses terminates with the cell or storm lifetime, whichever is the sooner. Associated with each pulse is a depth, $X$, so the pulse process is a marked point process (\cite{Cox80}).  The model developed by \cite{Cow07} allows pulse depths from a single cell to be dependent, but those from distinct cells are assumed independent.  No specific dependence structure is specified, and the model fitted in the paper assumed independent, exponentially distributed pulse depths, with mean depth $\mu_X$.  \cite{Cow07}'s model also assumed two superposed processes, with a common depth parameter across the two storm types, giving a total of eleven parameters.

\subsection{Development of the Random Parameter Bartlett-Lewis Instantaneous Pulse (BLIPR) model}\label{dev.random}
For the randomisation of $\eta$ in the BLIP model, we take the same approach as for the original Bartlett-Lewis model, but now with the additional assumption that the ratio, ($\omega = \xi/\eta$), of the pulse arrival rate to the cell duration parameter is kept constant.

In order to calculate the moments, it is helpful to think of the random parameter model as the superposition of a continuum of independent processes with random cell duration parameter, $\eta$, and storm origin rate, $\lambda f(\eta)$, where $f(\eta)$ is the density function of $\eta$.  Now, the $r$th cumulant of a sum of independent random variables is the sum of their $r$th cumulants.  Therefore the mean, variance and 3rd central moment (which are the first three cumulants) can simply be obtained by replacing $\lambda$ with $\lambda \,f(\eta)$ in their original equations, and integrating over possible values of $\eta$.

The integration approach described requires some expectations of functions of $\eta$.  In particular, we need $E_{\eta}\left[\left(\frac{1}{\eta}\right)^k e^{-\eta x}\right] $ for $k = 1$ and various values of $x$, given by:
\begin{align}
    E_{\eta}\left[\left(\frac{1}{\eta}\right)^k e^{-\eta x}\right] &= \frac{\nu^\alpha}{\Gamma(\alpha)}\int_0^{\infty} \eta^{\alpha-1-k} e^{-(\nu+x)\eta} d\eta \nonumber \\
    &= \frac{\nu^\alpha}{\Gamma(\alpha)} \times \frac{\Gamma(\alpha - k)}{(\nu+x)^{\alpha-k}}. \nonumber
\end{align}
Note that, in order for the integral not to diverge at zero, we require $\alpha > k$.  This proved to be an issue for the original Bartlett-Lewis model, as discussed in Section \ref{exist}, where the skewness integral included elements with $k = 4$.  For the BLIPR model, we only require $\alpha > 1$ in order for the integrals for the variance and skewness of the aggregated rainfall not to diverge.  However, the constraint $\alpha > 2$ may still be desirable (or a \ql truncated\qr \;version used) in order to prevent the simulation of unrealistically long rain events, as discussed.

The moments are derived from the original equations of \cite{Cow07}, by taking expectations over $\eta$ and using the formula above.  As in the original fixed parameter BLIP model, the flexibility to allow pulse depths to be dependent within cells is retained.  In their empirical fits, \cite{Cow07} assumed these to be independent, but intuitively, dependent pulse depths should allow higher values of extremes at short timescales.  This is desirable since the fits understated five-minute extreme values.  The moments for the new model are given in Appendix \ref{formulae}.

\section{Fitting methodology}\label{fitting}
The generalised method of moments (GMM) is used for fitting.  This is an extension of the method of moments which estimates parameters by equating expressions for population moments with their sample values. In the GMM, the number of properties to be fitted to exceeds the number of unknown parameters, and the estimator is given by the value of $\theta$ that minimises:
\[
    S(\theta|T) = (T - \tau(\theta))' W (T - \tau(\theta))
\]
for some positive definite weighting matrix $W$, where $\theta$ is the unknown parameter vector, $T$ is the vector of observed values for a set of $k$ properties, and $\tau(\theta)$ is the vector of their expected values under the model.  $S$ is referred to as the \ql objective function\qr.  \add[JK]{The optimal weights matrix (in terms of the identifiability of parameters) is the inverse of the covariance matrix of statistics \mbox{\citep{Han82}}, which here must be estimated empirically due to the complexity of the analytical expressions. In a recent simulation study using a point process based rainfall model, \mbox{\cite{Jesus11}} find that a two-step approach is required in order to derive a reliable sample estimate of the full covariance matrix, but that the diagonal matrix of inverse variances, calculated using just a single step, is close to optimal, and this is the approach that we follow.} \change[JK]{Here, we take $W$ to be a diagonal matrix, so that the}{The } objective function becomes $S(\theta|T) = \sum_{i=1}^k w_i [T_i(y) - \tau_i(\theta)]^2$, with the $w_i$ equal to $1/\text{Var}(T_i(y))$).  \add[JK]{Variances are calculated separately for each calendar month, pooling the data over observation years, and a separate fit is then produced for each month to allow for seasonality.}  \remove[JK]{This is a slight simplification of the theoretically optimal approach  (in terms of the identifiability of parameters) of \mbox{\cite{Han82}}, where $W$ is the inverse of the covariance matrix of statistics.}

Note that, since the number of properties included in $S$ exceeds the number of parameters, there is no guarantee that there will be a good fit to all the fitting properties.   The adequacy of the fit is thus assessed by considering properties used in the fitting procedure, as well as others that are of interest in hydrological applications. Some properties will need to be assessed using simulations, for example, extreme values.

We follow \cite{Cow07} in our choice of fitting properties - the hourly mean, plus the coefficient of variation, lag-1 correlation and skewness at timescales of 5 minutes, 1 hour, 6 hours and 24 hours.

Minimisation of $S$ requires a numerical optimisation routine. The approach followed here is that of \cite{Whe05}, and we have used the optimisation routines developed for that project.  Firstly, a set number of optimisations are carried out using the Nelder-Mead method, each starting with a different initial value for the set of parameters.  This set of initial values is generated by random perturbation about a single user-supplied value.  The best parameter set is then used as a new starting value for a further set of optimisations, which now use a Newton-type algorithm.  The reason for the use of two different optimisation routines is that the first is more robust and thus well suited to identifying promising regions of the parameter space, whereas the second is more powerful if given good starting values.

Different approaches to estimating parameter uncertainty have been taken in the literature. \cite{Rod88} look at parameter stability for the random parameter Bartlett-Lewis model by perturbing the input statistics by small amounts ($\pm$ 2\%) and looking at the impact on the resulting parameter estimates.  \cite{Cow98} uses a bootstrap approach, obtaining 100 sets of parameter estimates by fitting a Neyman-Scott model 100 times, each time using whole years sampled with replacement from the series of observed data.  Sampling whole years (separately for each calendar month) ensures that the dependencies in the rainfall series are captured.  \cite{Whe05} outline a method, based on the asymptotic theory of estimating equations, to estimate standard errors.  Another approach, used by \cite{Chan03}, is the examination of profile objective functions, which is our preferred method here.  Each parameter in turn is fixed at each of a set of values, and the objective function is optimised over the remaining parameters.  The resulting plot for each parameter showing the optimised objective function against the set of parameter values provides a useful means for assessing the identifiability of the parameter - for example, a very flat objective function indicates a wide range of plausible values.  \change[JK]{Approximate 95\% confidence intervals can also be calculated using the objective function itself (although there may be problems with numerical instabilities for the more complex models.)}{If the optimal weighting matrix is used, it can be shown that $2[S(\theta_p, \hat{\psi}|T) - S(\hat{\theta}_p, \hat{\psi}|T)]$ has a $\chi_1^2$ distribution, where $\theta_p$ is a single component of the parameter vector, $\theta$, and $\psi$ is the vector of the remaining parameters, and this result can be used to calculate approximate 95\% confidence intervals for each parameter.  Although this result does not hold if a non-optimal weighting matrix is used, a useful approximation has been suggested by \mbox{\cite{Jesus11}}}.

\section{Comparison of models on Bochum data}\label{comp}

\subsection{Models Fitted}\label{mods.fitted}
The models were fitted, using the methodology and fitting properties discussed, to 69 years of five-minute rainfall data from a single site in Bochum in Germany.  \remove[JK]{A separate fit was produced for each month, to allow for seasonality.}  In each case, we assume that $\sigma_X/\mu_X = 1$, and that E$[X^3] =  6 \mu_X^2$ (consistent with $X$ being exponentially distributed).  Initially, no further constraints were imposed on the parameters, other than that they should be greater than zero. The six models initially fitted were:

\textbf{Rectangular Pulse Models}
\begin{enumerate}
\item the Bartlett-Lewis Rectangular Pulse model (BLRP)
\item the Random Parameter Bartlett-Lewis Rectangular Pulse Model (BLRPR)
\item the Bartlett-Lewis Rectangular Pulse model with two superposed processes (BLRP2);
\end{enumerate}

\textbf{Instantaneous Pulse Models}
\begin{enumerate}
\item the Bartlett-Lewis Instantaneous Pulse Model (BLIP)
\item the Random Parameter Bartlett-Lewis Instantaneous Pulse model, introduced in Section \ref{dev.random} (BLIPR)
\item the Bartlett-Lewis Instantaneous Pulse model with two superposed processes (BLIP2).
\end{enumerate}

For the Bartlett-Lewis Rectangular Pulse model, on randomising the cell duration parameter, $\eta$, the fitted solution gave such a high precision to the mean cell duration, that it effectively replicated the non-random solution.  Thus, the fitted parameter set for the BLRPR model is simply a re-parameterised version of the set of BLRP parameters, and there is thus no improvement in the fit compared with the fixed $\eta$ version.  This appears to contradict examples in the literature where the randomised $\eta$ version had shown an improved fit compared to the fixed $\eta$ model \citep{Rod88, Whe05}.  On further investigation, we concluded that the improvement in the fit to proportion dry that had previously been found by randomising $\eta$ was at the expense of a deterioration in the fit to the skewness, which had not been included as a fitting property in these earlier analyses.  In particular, if skewness is not included in the fit, it is highly overestimated in the summer months at timescales of six and twenty-four hours.

Fitting the models with two superposed processes proved problematic.  Although the BLRP2 model with no parameter constraints gave a very good fit in terms of a low minimum objective function value, the parameters thus obtained were highly unstable, unrealistic and inconsistent from month to month, and no standard errors could be found.  It was clear that there was insufficient information in our observed data to identify the large number of required parameters.  Introducing constraints for the parameters increased the minimum objective function values, and did not resolve the situation, with resulting solutions having many parameters lying on the constraint boundaries. We therefore concluded that ensuring realistic and reasonably smooth parameters across months would require constraints on the relationships between parameters, rather than just setting bounds on individual parameters.  There were similar issues with the BLIP2 model.  Ultimately we decided that both of these models' parameter identifiability issues made them unsuitable for practical application\add[JK]{, at least in the context of this dataset and set of fitting properties, and our requirement for a  model that is both robust and  easy to fit}.

Given the above findings, we present further results here for the following three models only:  BLRP, BLIP, BLIPR.

For the BLIP and BLIPR models we initially followed \cite{Cow07} and assumed that pulses within a single cell had independent depths.  However, for the BLIPR model an alternative assumption was also considered, whereby pulses within a single cell have a common depth (the most extreme form of dependance).  The latter achieved a lower minimum objective function value in all months, and a better fit in respect of properties not included in the fitting process, such as wet/dry properties.  For both of these options,  the unconstrained solution gave an extremely high number of pulses per hour (of the order of $10^5$--$10^6$), so for practical reasons, $\mu_X$  was constrained to be $0.001$, reducing the number of parameters by one.  All other fitted parameters were broadly as before, except for a corresponding change in $\omega$.  The quality of the fit was unchanged with this constraint, as the product term $\mu_X \,\omega$ effectively forms a single composite parameter over most of the possible parameter space.  We also considered two alternative constraints on $\alpha$: $\alpha > 1$ or $\alpha > 2$, as discussed in Section \ref{dev.random}.  The former only affects July, whereas the latter affects all the summer months.

A comparison of the performance of the three fitted models, together with the findings discussed earlier in this section and  consideration of the fitted parameter sets (shown in Tables \ref{par.BLRP} to \ref{par.BLIPRd2} of Appendix \ref{app:params1})  led us to a hypothesis, which we present in the next section.

\subsection{Initial performance comparison of the fitted models}
Table \ref{min.obj} shows the minimum objective function value for each of the models that we have successfully fitted, for each month.  Since the same set of moments and weights were used for each model, these are directly comparable.
\begin{table}[ht]
\centering
\begin{small}
\begin{tabular}{lrcccccr} \toprule
&   &  BLIP$^1$   & BLIPR$^1$,  & BLIPR$^1$, & BLIPR$^2$,\\
&  BLRP &  independent &  independent  & common & common\\
& &  pulse depths  & pulse depths & pulse depths & pulse depths\\\midrule
Jan &	83	&	67		&	45 & 40 & 40\\
Feb	 &	38	&	56		&	30 & 24 & 24 \\
Mar	 &	100	&	113		&	58 & 48 & 48\\
Apr	 &	110	&	168		&	85 & 66 & 66\\
May	&	141	&	239		&	93 & 76 & 78\\
Jun	&	152	&	275		&	92 & 72 & 80\\
Jul	&	162	&	345		&	110 & 95 & 97\\
Aug	 &	140	&	268		&	86 & 76 & 79\\
Sep	 & 149	&	271		&	87 & 65 & 72\\
Oct	 & 92	&	150		&	71 & 50 & 50\\
Nov	 &	68	&	76		&	30 & 25 & 25\\
Dec	 &	68	&	67		&	32 & 28 & 28\\
\bottomrule
\end{tabular}
\end{small}
\caption[Comparison of minimum objective function values.]{Comparison of minimum objective function values; $^1$: $\alpha > 1$; $^2$: $\alpha > 2$.}\label{min.obj}.
\end{table}

Key findings from the results are summarised below:
\begin{itemize}
\item The BLRP model outperforms the BLIP model, with a lower minimum objective function value in all months except January and December.  The model with rectangular pulses has generally been considered unsuitable for timescales shorter than the mean cell duration, due to the unrealistic intensity shape.  However, when fine-scale data are available for fitting, the fitted model tends to have shorter, more frequent cells than if only hourly data are available (of the order of 5--10 minutes, compared with 20--40 minutes for most months), which are still within a realistic range.  With these shorter cells, and given also the potential for cells to overlap, repetition of the same rainfall totals over consecutive five minute intervals is relatively infrequent.  The fitted parameters are shown in Table \ref{par.BLRP} of Appendix \ref{app:params1}, along with some key properties such as mean storm and cell inter-arrival times and durations.
\item When skewness is included in the fit, there is no benefit to randomising the cell duration parameter in respect of the BLRP model, as discussed in Section \ref{mods.fitted}. However, there is a clear benefit in respect of the BLIP model, with the randomised version showing the best performance of all the models.
\item The fitted BLIPR model has a very high number of pulses per cell (particularly if we do not apply constraints, as discussed), with very short inter-arrival times, and the better performing version has common within-cell pulse depths.  Effectively then, the cells are \ql rectangular\qr.
\end{itemize}
These results imply that it is not the replacement of rectangular pulses by clusters of instantaneous ones that leads to the improved performance of the BLIPR model, compared with the BLRPR model. Instead, the improved performance can be attributed to the fact that the BLIPR model allows rainfall intensity to vary with cell duration, since the pulse rate effectively drives the intensity and is proportional to the cell duration parameter, $\eta$. Our new model variant thus gives a simple, but effective way of introducing dependence between cell duration and intensity.

This suggests that the same effect could be achieved by amending the BLRPR model, so that the mean cell intensity parameter, $\mu_X$ is also varied in proportion to the cell duration parameter, $\eta$.  This is preferable from a computational point of view, eliminating the need for simulation of a vast number of instantaneous pulses.

\subsection{Testing our hypothesis and new model variant}\label{test.hyp}
Extending the BLRPR model to allow $\mu_X$ to vary in proportion to the cell duration parameter, $\eta$, is  straightforward, and  follows the methodology discussed in Section \ref{dev.random}.  We re-parameterise the BLRPR model so that the ratio, $\iota = \mu_X/\eta$ is now kept constant, and  express $E(X^2)$ and $E(X^3)$ in terms of $\iota$ also (for which the formulae depend on the choice of distribution for the rainfall intensity).  We then take expectations over $\eta$ as before. The analytical expressions for this new model, which we denote the BLRPR$_X$  model, are given in Appendix \ref{formulae}.

The fitted parameter set, assuming an exponential distribution for cell intensities as before, is given in Table \ref{par.BLRPRM2} of Appendix \ref{app:params1}.  Comparing this with the fitted parameters of the Random Parameter Bartlett-Lewis Instantaneous Pulse (BLIPR) model with common within-cell pulse depths, shown in Table \ref{par.BLIPRd2}, the strong similarity between the two models is evident.  In particular, the new parameter, $\iota$, of the $\text{BLRPR}_X$ model broadly equates to $\mu_X \, \omega$ of the BLIPR model (noting that $\mu_X$ represents an intensity in the rectangular pulse models, but a depth in those with instantaneous pulses).  Values of the minimum objective function (see Table \ref{min.obj2})  and plots of the two fits are also found to match, thus supporting our hypothesis.
\begin{table}[ht]
\centering
\begin{small}
\begin{tabular}{lcccccccccccc} \toprule
&Jan & Feb & Mar & Apr & May & Jun & Jul & Aug & Sep & Oct & Nov & Dec\\ \midrule
BLIPR & 40 &	24 &	48 &	66 &	79 &	80&	97&	79	&72	&50&	25&	 28\\
$\text{BLRPR}_X$ &39	&22	&46	&63	&74	&76	&92	&74	&68	&47	&23	&26\\
\bottomrule
\end{tabular}
\end{small}
\caption[Comparison of minimum objective function value; $\alpha > 2$.]{Comparison of minimum objective function value; $\alpha$ constrained to be at least 2.}\label{min.obj2}
\end{table}

We have therefore established that the new rectangular pulse model variant is effectively equivalent to the BLIPR model with common within-cell pulse depths, and that there is therefore no need to replace the rectangular pulses with a process of  instantaneous pulses for fine-scale data.  This is the optimal model, at least in terms of the minimum objective function values.   In the next section we examine the performance of the three models (BLRP, BLIP, $\text{BLRPR}_X$) in more detail, firstly in terms of the fitted moments, and then by considering wet/dry properties, which were not included within the objective function, and extreme value performance.

\subsection{Performance comparison of the fitted models}

\subsubsection{Fitted Moments\label{fit.moments}}
Plots of the fits of the models (BLRP, BLIP, $\text{BLRPR}_X$) against the observed data for each month in respect of the mean, coefficient of variation, lag-1 autocorrelation and skewness coefficient are shown in Figures \ref{mean}-\ref{skew} in Appendix \ref{plots}.   Note that the $y$-axes for these and other similar plots in this paper have been selected automatically such that, for each individual plot,  the axis spans the range covered by the observed and fitted values.  This means that the fit in respect of an individual model tends to look worse if all models fit well, than if at least one of the other plotted models has a poor fit (since in the latter case the scale will be wider).  Care should therefore be taken to consider also the scale when examining such plots.

All the models generally perform well with respect to the properties included in the fitting. They reproduce the mean exactly (this is not a given, since the number of properties fitted exceeds the number of parameters), and fit the coefficient of variation well at all timescales.  All tend to underestimate the lag-1 autocorrelation at longer timescales.  All also tend to underestimate the skewness at the shorter timescales, with the $\text{BLRPR}_X$ model showing the best fit in respect of 5 minute skewness, and the BLIP model the worst.

\subsubsection{Wet/dry properties\label{fit.wet.dry}}
The proportion of dry intervals is a very important property for hydrological applications.  Although this could have been included as one of the fitting properties, it is useful to reserve an important feature for subsequent model validation, as this gives an independent  test of the appropriateness of the model structure.  Plots of the fits of the models against the observed data for each month in respect of the proportion dry are shown in Figure \ref{pdry}.  The $\text{BLRPR}_X$ model can be seen to outperform the other models with respect to the fit to proportion dry, across all timescales.

It is also of interest to consider the wet and dry spell transition probabilities (i.e the probability that a wet interval is followed by another wet interval, or a dry by another dry), which are important for the accurate modelling of antecedent conditions.  Figure \ref{ww} shows that the $\text{BLRPR}_X$ model again outperforms the other models with respect to the wet spell transition probability.  While the BLRP model has a good fit at the hourly timescale, it performs poorly at other timescales, with only the $\text{BLRPR}_X$ model showing consistency of performance across timescales.   There is less difference between models for the dry spell transition probabilities, with all models providing a reasonable fit at all timescales.  The fit to the wet/dry properties in respect of the summer months would be further improved if we did not impose the constraint that $\alpha > 2$.  However, this would be at the expense of allowing storms and cells of unrealistic durations in the simulations (as discussed in Section \ref{dev.random}), and also a slight deterioration of the fit to the 24 hour variance, and 6 hour lag-1 autocorrelation.

\subsubsection{Extreme value performance\label{fit.extremes}}
For our data, the months with the highest rainfall, rainfall variability and skewness are the summer months, and these are also the months with the highest extremes.  A comparison of the fit of extremes for July for the $\text{BLRPR}_X$ model is given in Figure \ref{ext.jul}, using Gumbel plots.  We use 100 simulations, allowing for parameter uncertainty by sampling from the distribution of the parameter set for each simulation, which is of the same length as the observed data.  The maximum rainfall per unit-time is plotted against the \ql reduced-variate\qr \, $-\ln(-\ln(1-1/R))$ where $R$ is the return period. The graphs for July show that the model has a tendency to underestimate extremes, as has been noted before for this type of model.  Results for other months give a fairly similar picture.

A comparison showing mean annual extremes (averaged over fifty simulations) for a number of alternative models at the five minute and hourly timescales is also shown in Figure \ref{ext.mean}.  At the five minute timescale, the $\text{BLRPR}_X$ model gives the best performance, although all the models underestimate the extremes.  Results are closer at the one-hour timescale, and for longer timescales, there is essentially no difference between models.

Based on our analysis, the $\text{BLRPR}_X$ is shown to be the best performing of the models compared, both in terms of the moments fitted, and more importantly, in respect of the wet/dry properties and extreme values, neither of which is included in the fit. It is also intuitively appealing, since the intensity of rainfall is known to vary inversely with the duration of the rain event.  Further, this dependence has been introduced to the BLRPR model without the need for any additional parameters or complexity.  Considering the fitted parameter set, shown in Table \ref{par.BLRPRM2} of Appendix \ref{app:params1}, the parameter values change fairly smoothly from month to month.  Comparing with empirical observations from \cite{Hou82}, the parameter values seem reasonable.  Winter storms last several hours, have around 20 cells, which each last on average around 22 minutes.  In summer, storms and cells are shorter, and have around 8 cells. However, these have a correspondingly much higher intensity, giving broadly the same amount of rainfall per storm over all months.

\section{Parameter Identifiability and Confidence Intervals}\label{param.id}
Finally, we explore the parameter identifiability of the new $\text{BLRPR}_X$ model using profile objective functions as described in Section \ref{fitting}.  Profile objective functions for the logarithm of the parameters are shown in Figure \ref{profiles} for the month of January.  As before, we have constrained the value of $\alpha$ to be greater than 2 (except in the plot for $\alpha$ itself).  The first set of plots shows a wide range of possible parameter values, allowing us to check whether there are multiple local minima, for example, or extensive regions where the objective function is flat.  We have then reduced the parameter range so that the approximate  95\% confidence intervals can be seen more clearly.  The plots show that the parameters are fairly well identified in January.  Results for other months (not shown) again indicate good parameter identification, although there is slightly more parameter uncertainty in the summer months, which is fairly typical.

\section{Potential further improvements}\label{improve}
Although the new model fits well, there are areas for improvement, the most important of which is the fit to extreme values, which are understated.  This is perhaps surprising, as intuitively the inclusion of the skewness coefficient as one of the fitting properties should lead to an improved fit in respect of extremes.  On investigation, we found that our approach of averaging the skewness over 69 separate observation months, rather than calculating a single statistic over the whole of the data, tends to understate the skewness coefficient itself, particularly at the 5 minute timescale.  This is found to be related to the effective weights that are applied to periods of high skewness under the two alternative approaches, rather than to sampling variation, or to the choice of mean (local or global) about which the moments are centred.  The alternative approach would have given us a slightly better fit to extremes, but does not permit the calculation of the covariance matrix for the observed statistics, since it is just a single sample.

It is generally thought that a distribution for the rainfall intensity with \ql fatter tails\qr \; should improve the fit to extremes. This was investigated by replacing the exponential distribution with the gamma and Weibull distributions, both of which have the exponential as a special case.  The Weibull distribution gave the better performance here.  However, the addition of a further parameter caused problems in terms of parameter identifiability, with less consistency from month to month.  Also the asymptotic results showed a very high correlation between the estimated shape parameter and the intensity parameter, suggesting that an additional parameter is not justified.  Constraining the shape parameter to a fixed value may, however, be a viable strategy, since this does not require any additional parameters.  Here the fitted shape parameter was close to 0.6 in most months.  Using a Weibull conditional intensity distribution with a shape parameter of 0.6 rather than an exponential (which has a  shape parameter of 1) improves the fits to 5 minute skewness and to extremes at short timescales, but with some deterioration in the lag-1 autocorrelation at longer timescales.

Another alternative considered for the $\text{BLRPR}_X$ model involved allowing a more flexible intensity/duration relationship, by letting the mean intensity be proportional to the cell duration parameter, raised to some fixed power, the level of which is to be determined i.e. to have $\iota = \mu_X/\eta^c$ for some additional parameter, $c$.  However, it was found that the fitted values of $c$ were fairly close to 1, suggesting that this additional complexity is not required.  Alternative formulae for this relationship could be considered, but are likely to affect the analytical tractability.

\section{Discussion and Conclusions}\label{conclusions}
In this paper, using a structured comparison of different versions of the Bartlett Lewis clustered point process-based models, including two new model extensions, we have clarified some key aspects of performance.  Our focus here has been on fine scale data.  We have highlighted some limitations in all the models, notably an inability to achieve a good fit to all properties in the summer months of a temperate climate, when rainfall exhibits particularly high variability and skewness.  Such limitations are not surprising when we consider the simplicity of these models compared with the highly complex real physical rainfall process.  The challenge of achieving a good fit at timescales that cover the wide range from five minutes up to daily, is particularly demanding.   However, the performance of the original rectangular pulse (BLRP) model, originally considered unsuitable for fine-scale data due to the unrealistic rectangular pulses, has far exceeded prior expectations.

We showed that the main driver behind improved performance, particularly in respect of skewness and extremes at short timescales, is the introduction of an inverse dependence between rainfall intensity and cell duration. Our proposed new model, which is an extension of the Random Parameter Bartlett-Lewis Rectangular Pulse model, gives a simple but effective way of introducing such dependence, with no increase in the number of parameters.  It allows the rainfall intensity parameter, $\mu_X$, previously assumed to be constant,  to vary in proportion to the cell duration parameter, $\eta$, which itself varies randomly between storms.   Although  instantaneous pulses were useful in leading us to this conclusion, ultimately we discovered that they are not required, and the computationally simpler rectangular pulse version is preferred.  Adding further parameters adds little to the fit, since typically improvements in some properties cause degradation in others, and more parameters bring issues of parameter identifiability and consistency.  Replacing the exponential intensity distribution with a Weibull with a fixed shape parameter, however, may be desirable.  \add[JK]{The introduction of depth-duration dependence is desirable for all datasets and timescales, although the positive impact of the new model is greatest for fine-scale data}.

Appropriately allowing for the uncertainty in the parameter estimates is important, and this is often not addressed in the hydrological literature.  A suggested approach, used here when generating simulations, is to sample parameters from the asymptotic multivariate normal distribution (or lognormal if the logarithms of the parameters are fitted).  In this way, rare, but potentially damaging scenarios should be better represented in simulations, particularly if, as here, extreme values tend to be underestimated by the model.

\section{Acknowledgements}
Deutsche Montan Technologie and Emschergenossenschaft/Lippeverband in Germany are gratefully acknowledged for providing the data.  We would also like to thank Richard Chandler and Joao Jesus for helpful advice.  Jo Kaczmarska is also pleased to acknowledge financial support from the Engineering and Physical Sciences Research Council.

\pagebreak
\begin{appendix}
\appendixpage
\section{Formulae for fitting properties}\label{formulae}
In this Appendix, we give the formulae for the mean, variance, lag-1 autocovariance  and 3rd central  moment of the discrete-time aggregated process in respect of the BLIPR and $\text{BLRPR}_X$ models. Throughout, the timescale to which the continuous process is aggregated is denoted as $h$. The required fitting properties can be derived from those given here, as follows:
\begin{small}
\begin{align}
    \text{Coefficient of variation} = \frac{\sqrt{\text{Var}[Y_i^h]}}{\text{E}[Y_i^h]},
\end{align}
\begin{align}
    \text{Skewness coefficient} = \frac{\text{E}[ (Y_i^h- \text{E}(Y_i^h))^3 ] }{\text{Var}[Y_i^h]^{3/2}},
\end{align}
\begin{align}
    \text{Lag 1 autocorrelation} = \frac{\text{Cov}(Y_i^h, Y_{i+1}^h)}{\text{Var}[Y_i^h]}.
\end{align}
\end{small}

The random cell intensity, denoted $X$, has been assumed to have a one parameter distribution.  The parameterisation in respect of the Exponential distribution, used here, is as follows:

\begin{tabular}{llll}
\textit{Parameter}:  &   $\mu_X$ &&\\
\textit{Moments}:   &     E$(X) = \mu_X$ &  E$(X^2) =  2 \mu_X^2$ & E$(X^3) = 6 \mu_X^3$\\
\end{tabular}

\subsection{Moments for the Barlett-Lewis Instantaneous Pulse Random $\eta$ (BLIPR) model}\label{A.BLIPR}
\textbf{Parameter definitions}\vspace{-10pt}
\begin{small}
\begin{itemize}
\item $\lambda$ - storm arrival rate
\item $\alpha$ - shape parameter for the Gamma distribution of the cell duration parameter, $\eta$
\item $\nu$ - scale parameter for the Gamma distribution of $\eta$
\item $\kappa$ - ratio of the cell arrival rate to $\eta$ (i.e. $\beta/\eta$)
\item $\phi$ - ratio of the storm (cell process) termination rate to $\eta$ (i.e. $\gamma/\eta$)
\item $\omega$ -  ratio of the pulse arrival rate to $\eta$ (i.e. $\xi/\eta$)
\item $\mu_X$ - mean pulse depth
\item E$(X^2)$ - mean of squares of pulse depths
\item E$(X^3)$ - mean of cubes of pulse depths
\item E$(X_{ijk}X_{ijl})$ - product moment of the depths of 2 pulses within the same cell
\item E$(X_{ijk}X_{ijl}X_{ijm})$ - product moment of the depths of 3 pulses within the same cell
\item $\mu_p = \frac{\kappa \omega}{\phi(\phi + 1)}$ - mean number of pulses per storm
\end{itemize}

\textbf{Mean}\vspace{-10pt}
\begin{align}
\text{E}[Y_i^h] &= \lambda \mu_p \mu_X h.
\end{align}
\textbf{Variance}
\begin{align}
\text{Var}[Y_i^h] &= \lambda \mu_p  \Bigg\{ \text{E}(X^2)h + \frac{2\mu_X^2 \kappa \omega}{\phi^2} \text{E}_{\eta}\Bigg(\frac{1}{\eta} e^{-\phi \eta h}- \frac{1}{\eta} + \phi h\Big)  \nonumber \\
 &   + \frac{2 \omega}{(\phi+1)^2} \Bigg[ \text{E}(X_{ijk}X_{ijl}) - \mu_X^2 \kappa \frac{\phi}{\phi + 2} \Bigg] \text{E}_{\eta}\Bigg(\frac{1}{\eta} e^{-(\phi+1)\eta h} - \frac{1}{\eta} + (\phi + 1)h\Bigg) \Bigg\}\nonumber\\\nonumber \\
&= \lambda \mu_p  \Bigg\{ \text{E}(X^2)h + \frac{2\mu_X^2 \kappa \omega}{\phi^2} \Bigg(\frac{\nu^\alpha}{(\alpha-1)(\nu + \phi h)^{\alpha-1}}- \frac{\nu}{\alpha-1} + \phi h\Bigg)  \nonumber \\
&  + \frac{2 \omega}{(\phi+1)^2} \Bigg[ \text{E}(X_{ijk}X_{ijl}) - \mu_X^2 \kappa \frac{\phi}{\phi + 2} \Bigg] \Bigg(\frac{\nu^\alpha}{(\alpha-1)(\nu+(\phi+1) h)^{\alpha-1}} - \frac{\nu}{\alpha-1} \nonumber \\
& + (\phi + 1)h\Bigg) \Bigg\}.
\end{align}
\textbf{Covariance  at lag $k \geq 1$}
\begin{align}
\lefteqn{\text{Cov}(Y_i^h, Y_{i+k}^h)  } \nonumber \\
 &= \lambda \mu_p \omega \Bigg[ \frac{\mu_X^2 \kappa}{\phi^2} \text{E}_{\eta}\Bigg(\frac{e^{-\phi \eta (k-1) h} - 2 e^{-\phi \eta k h} + e^{-\phi \eta (k+1) h}}{\eta}\Bigg) \nonumber \\
 & \;\;+ \Bigg( \text{E}(X_{ijk} X_{ijl}) - \mu_X^2 \kappa \frac{\phi}{(\phi + 2)}\Bigg)\text{E}_{\eta}\Bigg(\frac{e^{-(\phi+1)\eta (k-1)h}  - 2 e^{-(\phi+1)\eta kh} + e^{-(\phi+1)\eta (k+1) h}}{(1+\phi)^2 \eta}\Bigg)\Bigg] \nonumber \\\nonumber \\
&= \lambda \mu_p \omega \Bigg(\frac{\nu}{\alpha - 1}\Bigg)\Bigg[ \frac{\mu_X^2 \kappa}{\phi^2} \Bigg\{\Bigg(\frac{\nu}{\nu+\phi (k-1) h}\Bigg)^{\alpha-1}  \nonumber \\
& \;\;\;- 2\Bigg(\frac{\nu}{\nu + \phi k h}\Bigg)^{\alpha-1} + \Bigg(\frac{\nu}{\nu+\phi (k+1) h}\Bigg)^{\alpha-1}\Bigg\} + \Bigg( \text{E}(X_{ijk} X_{ijl}) - \mu_X^2 \kappa \frac{\phi}{(\phi + 2)}\Bigg) \nonumber \\
 &\;\;\; \times \Bigg\{\Bigg(\frac{\nu}{\nu+(\phi+1) (k-1) h}\Bigg)^{\alpha-1} - 2\Bigg(\frac{\nu}{\nu + (\phi+1) k h}\Bigg)^{\alpha-1} \nonumber \\
 &\;\;+ \Bigg(\frac{\nu}{\nu+(\phi+1) (k+1) h}\Bigg)^{\alpha-1}\Bigg\} \Bigg].
\end{align}
\textbf{3rd central moment}
\begin{align}
\lefteqn{\text{E}[ (Y_i^h- \text{E}(Y_i^h))^3 ]  } \nonumber \\
&=  \lambda \kappa  \omega^3 \Bigg\{ \frac{6}{(1 + \phi)^3} \Bigg[\frac{\text{E}(X_{ijk}X_{ijl}X_{ijm})}{\phi } + \frac{2 \text{E}(X_{ijk}X_{ijl}) \mu_X \kappa}{\phi (2+\phi)}  - \frac{\mu_X^3 \kappa^2}{(2+\phi)}\Bigg]\nonumber \\
 &\times \Bigg[h - \frac{2\nu}{(\alpha-1)(1+\phi)} + \frac{2 \nu }{(1+\phi)(\alpha-1)}\Bigg(\frac{\nu}{\nu+(1+\phi)h}\Bigg)^{\alpha-1} + h \Bigg(\frac{\nu}{\nu+(1+\phi)h}\Bigg)^\alpha\Bigg] \nonumber \\
& + \frac{6}{(1+\phi)(2 + \phi)^2}\Bigg[-\frac{2 \text{E}(X_{ijk}X_{ijl}) \mu_X \kappa}{ (1+\phi)} + \frac{\mu_X^3 \kappa^2}{(3+\phi)}\Bigg]\nonumber \\
&  \times \Bigg[h-\frac{\nu}{(\alpha-1)}\Bigg\{\frac{3+2 \phi}{(1+\phi)(2+\phi)} - \Bigg(\frac{2+\phi}{1+\phi}\Bigg)\Bigg(\frac{\nu}{\nu + (1+\phi)h}\Bigg)^{\alpha-1}  \nonumber \\
&+ \Bigg(\frac{1+\phi}{2 + \phi}\Bigg)\Bigg(\frac{\nu}{\nu+(2+\phi) h}\Bigg)^{\alpha-1}\Bigg\}\Bigg] \nonumber \\
&+ \frac{6 \mu_X^3 \kappa^2}{ \phi^3 (1+\phi)}\Bigg[h - \frac{2 \nu}{\phi (\alpha-1)} + \frac{2 \nu}{\phi (\alpha - 1)}\Bigg(\frac{\nu}{\nu + \phi h}\Bigg)^{\alpha-1} + h \Bigg(\frac{\nu}{\nu + \phi h}\Bigg)^\alpha\Bigg] \nonumber \\
&+  \frac{6}{\phi(1+\phi)^2} \Bigg[\frac{2 \text{E}(X_{ijk}X_{ijl}) \mu_X \kappa}{ \phi } - \frac{\mu_X^3 \kappa^2}{ (2+\phi)}\Bigg] \nonumber \\
& \times \Bigg[h-\frac{\nu}{(\alpha-1)}\Bigg\{\frac{1+2\phi}{\phi (1+\phi)} - \frac{(1+\phi)}{\phi}\Bigg(\frac{\nu}{\nu+ \phi h}\Bigg)^{\alpha-1} + \frac{\phi}{(1+\phi)}\Bigg(\frac{\nu}{\nu + (1+\phi)h}\Bigg)^{\alpha-1}\Bigg\}\Bigg] \nonumber \\
&+ \frac{6 \text{E}(X_{ijk}^2 X_{ijl})}{ \omega \phi (1+\phi)^2}\Bigg[h - \frac{\nu}{(1+\phi)(\alpha-1)}\Bigg\{ 1 - \Bigg( \frac{\nu}{\nu + (1+\phi)h}\Bigg)^{\alpha-1}\Bigg\}\Bigg] \nonumber \\
&+  \frac{6 \text{E}(X^2) \mu_X \kappa}{\omega \phi^2 (1 + \phi)} \Bigg[h - \frac{\nu}{\phi (\alpha-1)} + \frac{\nu}{\phi (\alpha-1)}\Bigg(\frac{\nu}{\nu+\phi h}\Bigg)^{\alpha-1} - \frac{\phi^2}{(1+\phi)(2+\phi)}\nonumber \\
&\times \Bigg( h - \frac{\nu}{(1+\phi)(\alpha-1)} + \frac{\nu}{(1+\phi)(\alpha-1)}\Bigg(\frac{\nu}{\nu + (1+\phi)h}\Bigg)^{\alpha-1} \Bigg)\Bigg] + \frac{\text{E}(X^3)h}{\omega^2 \phi (1+\phi)}   \Bigg\}.\nonumber\\
\end{align}

\subsection{Moments for the Random Parameter Bartlett-Lewis Rectangular Pulse model with dependent intensity-duration ($\text{BLRPR}_X$)}\label{A.BLRPRM}
All expectations are left in the form $E_{\eta}\left[\eta^{-k}\;e^{-\eta s}\right] $ for various values of $k$ and $s$, and may be evaluated as:
\begin{align}
    E_{\eta}\left[\eta^{-k} e^{-\eta s}\right]
    &= \frac{\nu^\alpha}{\Gamma(\alpha)} \times \frac{\Gamma(\alpha - k)}{(\nu+s)^{\alpha-k}}, \;\;\;\;\;\text{for } \alpha > k.\nonumber
\end{align}
\textbf{Parameter definitions}\vspace{-10pt}
\begin{itemize}
\item $\lambda$ - storm arrival rate
\item $\alpha$ - shape parameter for the Gamma distribution of the cell duration parameter, $\eta$
\item $\nu$ - scale parameter for the Gamma distribution of $\eta$
\item $\kappa$ - ratio of the cell arrival rate to $\eta$ (i.e. $\beta/\eta$)
\item $\phi$ - ratio of the storm (cell process) termination rate to $\eta$ (i.e. $\gamma/\eta$)
\item $\iota$ - ratio of mean cell intensity to $\eta$ (i.e. $\mu_X/\eta$)
\item $f_1$ - $\text{E}(X^2)/\mu_x^2$
\item $f_2$ - $\text{E}(X^3)/\mu_x^3$
\item $\mu_C = 1 + \kappa/\phi$ - mean number of cells per storm
\end{itemize}

\textbf{Mean}\vspace{-10pt}
\begin{align}
\text{E}[Y_i^h] = \lambda h \iota \mu_c.
\end{align}
\textbf{Variance}
\begin{align}
\text{Var}[Y_i^h] &= 2 \lambda \mu_c \iota^2 \Bigg\{ \Bigg\{f_1 + \frac{\kappa}{\phi} \Bigg\} h + \text{E}_\eta (\eta^{-1})\Bigg\{\frac{\kappa (1-\phi^3)}{\phi^2 (\phi^2 - 1) } - f_1 \Bigg\} \nonumber \\ &\;\;- {} \text{E}_\eta (\eta^{-1}  e^{-\phi \eta h}) \frac{\kappa}{\phi^2 (\phi^2 - 1)} + \text{E}_\eta (\eta^{-1} e^{-\eta h}) \Bigg\{f_1 + \frac{\kappa \phi }{\phi^2 - 1}\Bigg\}\Bigg\}.
\end{align}
\textbf{Covariance  at lag $k \geq 1$}
\begin{align}
\text{Cov}(Y_i^h, Y_{i+k}^h)  &= \lambda \mu_c \iota^2 \Bigg\{ \Bigg(f_1 + \frac{\kappa \phi}{\phi^2 - 1}\Bigg)\bigg[ \text{E}_\eta (\eta^{-1} e^{-\eta (k-1) h}) - 2 \,\text{E}_\eta(\eta^{-1} e^{-\eta k h}) \nonumber \\
& \;\; + {}\text{E}_\eta(\eta^{-1} e^{-\eta (k+1) h}) \bigg]
 - \frac{\kappa} {\phi^2 (\phi^2 - 1)} \Bigg[\text{E}_\eta(\eta^{-1} e^{-\phi \eta (k-1) h}) - 2 \, \text{E}_\eta(\eta^{-1} e^{-\phi \eta k h}) \nonumber \\
 &\;\; +  \text{E}_\eta(\eta^{-1} e^{-\phi \eta (k+1) h})\Bigg]\Bigg\}.
\end{align}
\textbf{3rd central moment}
\begin{align}
\lefteqn{\text{E}_\eta[ (Y_i^h- \text{E}(Y_i^h))^3 ]  } \nonumber \\
&= \frac{\lambda \mu_c \iota^3 }{(1+2\phi+\phi^2)
    (\phi^4-2\phi^3-3\phi^2+8\phi-4) \phi^3} \nonumber \\
&\times  \Bigg\{  \text{E}_\eta\left[\eta^{-1} e^{-\eta h}\right] \Bigg(12 \phi^7  \kappa^2 - 24 f_1 \phi^2 \kappa - 18 \phi^4 \kappa^2 + 24 f_1 \phi^3 \kappa
- 132 f_1 \phi^6\kappa   + 150f_1 \phi^4 \kappa \nonumber \\
& - 42\phi^5 \kappa^2 - 6 f_1 \phi^5  \kappa + 108\phi^5 f_2
- 72 \phi^7 f_2  - 48 \phi^3 f_2 + 24 f_1 \mu_x \phi^8 \kappa   + 12\phi^3  \kappa^2 + 12\phi^9 f_2\Bigg) \nonumber \\
&+  \text{E}_\eta\left[ e^{-\eta h}\right] \Bigg(24 f_1 \phi^4 h \kappa + 6 \phi^9 h f_2- 30 f_1 \phi^6 h \kappa + 6 f_1 \phi^8 h \kappa + 54 \phi^5 h f_2 - 24 h f_2 \phi^3 - 36 \phi^7 h f_2\Bigg) \nonumber \\
&+  \text{E}_\eta\left[\eta^{-1} e^{-\eta \phi h}\right] \Bigg(- 48 \kappa^2 + 6 f_1 \phi^4 \kappa  - 48 \phi f_1 \kappa + 6\phi^5 \kappa^2 - 24f_1 \phi^2 \kappa + 36f_1 \phi^3  \kappa\nonumber \\
 & - 6f_1 \phi^5 \kappa  +  84\phi^2  \kappa^2 + 12\phi^3  \kappa^2 - 18\phi^4 \kappa^2 \Bigg)\nonumber  \\
 &+  \text{E}_\eta\left[ e^{-\eta \phi h}\right] \Bigg(- 24\phi h \kappa^2 +
 30\phi^3 h\kappa^2 - 6\phi^5 h  \kappa^2 \Bigg)\nonumber \\
 &+  \text{E}_\eta\left[\eta^{-1} \right] \Bigg(72\phi^7 f_2 + 48\phi f_1 \kappa +
 24 f_1 \phi^2 \kappa - 36f_1 \phi^3 \kappa - 84 \phi^2 \kappa^2 + 6 f_1 \phi^5 \kappa + 117 f_1\phi^6 \kappa   \nonumber \\
 &  + 39 \phi^5 \kappa^2  - 12 \phi^9 f_2  - 138 f_1 \phi^4 \kappa + 48 \kappa^2 - 9 \phi^7  \kappa^2 + 48\phi^3 f_2  + 18\phi^4 \kappa^2  - 21\phi^8 f_1 \kappa \nonumber \\
  & - 12 \phi^3  \kappa^2 - 108 \phi^5 f_2\Bigg)\nonumber \\
 &+  \Bigg(- 24 \phi h  \kappa^2  - 72f_1 \phi^6 h \kappa   - 36 \phi^5 h  \kappa^2 + 54 \phi^3 h \kappa^2 + 6 \phi^7 h \kappa^2 + 54\phi^5 h f_2 - 36 \phi^7 h f_2  \nonumber \\
 &- 24 \phi^3 h f_2  - 48 f_1 \phi^2 h \kappa  + 12f_1\phi^8 h \kappa + 6\phi^9 h f_2 + 108 f_1 \phi^4 h \kappa \Bigg) \nonumber \\
 &+ \text{E}_\eta\left[\eta^{-1} e^{-2 \eta h}\right] \Bigg(- 12 f_1\phi^4 \kappa - 3f_1 \phi^8 \kappa + 15 f_1 \phi^6 \kappa  - 3\phi^7  \kappa^2  + 3\phi^5  \kappa^2\Bigg) \nonumber\\
 &+ \text{E}_\eta\left[\eta^{-1} e^{-\eta h (1+\phi)}\right] \Bigg(- 24 f_1 \kappa \phi^3  - 6 f_1\phi^4 \kappa  + 6\phi^5 f_1 \kappa + 24 f_1 \kappa \phi^2  + 18 \phi^4 \kappa^2 \nonumber \\
 &- 12 \phi^3  \kappa^2 - 6 \phi^5\kappa^2 \Bigg)   \Bigg\}.
\end{align}
\end{small}
\clearpage
\section{Fitted parameters}\label{app:params1}
For each model, as well as the fitted parameters, we show a number of key properties, in order to allow a better comparison of models with different parameterisations.  The acronyms used for these properties are given below:

\begin{small}
\begin{tabular}{ll}
MSIT & mean storm inter-arrival time, hours\\
MSD & mean duration of storm activity, hours\\
MCIT & mean cell inter-arrival time, minutes\\
MCD & mean cell duration, minutes \\
MCS & mean number of cells per storm ($=\mu_C$) \\
MPC & mean number of pulses per cell \\
\end{tabular}
\end{small}

\begin{table}[ht]
\begin{footnotesize}
\begin{tabular}{lrrrrrrrrrrr} \toprule
	&	$\lambda$	&	$\mu_X$&	$\beta$	&	$\gamma$	&	$\eta$	&&	 MSIT  &		MSD &		MCIT &	MCD	&	MCS		\\ \midrule
Jan	&	0.022	&	0.960	&	5.422	&	0.231	&	5.975	&&	45.0	&	 4.3	&	11.1	&	10.0	&	24.5\\
Feb	&	0.021	&	0.942	&	5.142	&	0.260	&	5.310	&&	47.1	&	 3.8	&	11.7	&	11.3	&	20.7\\
Mar	&	0.021	&	1.334	&	4.478	&	0.262	&	7.061	&&	47.2	&	 3.8	&	13.4	&	8.5	&	18.1\\
Apr	&	0.022	&	1.944	&	3.829	&	0.271	&	8.387	&&	45.7	&	 3.7	&	15.7	&	7.2	&	15.1\\
May	&	0.023	&	3.662	&	3.157	&	0.370	&	9.239	&&	44.3	&	 2.7	&	19.0	&	6.5	&	9.5\\
Jun	&	0.025	&	6.431	&	2.694	&	0.413	&	11.154	&&	39.2	&	 2.4	&	22.3	&	5.4	&	7.5\\
Jul	&	0.023	&	10.136	&	1.672	&	0.356	&	12.011	&&	43.5	&	 2.8	&	35.9	&	5.0	&	5.7\\
Aug	&	0.023	&	7.072	&	2.411	&	0.408	&	11.066	&&	43.4	&	 2.5	&	24.9	&	5.4	&	6.9\\
Sep	&	0.021	&	5.306	&	2.945	&	0.379	&	10.470	&&	47.1	&	 2.6	&	20.4	&	5.7	&	8.8\\
Oct	&	0.019	&	2.209	&	4.071	&	0.275	&	8.104	&&	53.3	&	 3.6	&	14.7	&	7.4	&	15.8\\
Nov	&	0.023	&	1.207	&	5.884	&	0.276	&	6.741	&&	42.8	&	 3.6	&	10.2	&	8.9	&	22.3\\
Dec	&	0.024	&	1.059	&	5.475	&	0.265	&	5.906	&&	41.1	&	 3.8	&	11.0	&	10.2	&	21.7\\
\bottomrule
\end{tabular}
\caption{Parameters for BLRP model.}\label{par.BLRP}
\end{footnotesize}
\end{table}

\begin{table}[ht]
\begin{footnotesize}
\begin{tabular}{lrrrrrrrrrrrrrr} \toprule
	&	$\lambda$	&	$\mu_X$	&	$\beta$	&	$\gamma$	&	$\eta$	&	 $\xi$	&&	MSIT	&	MSD	&	MCIT	&	MCD	&	MCS	& 		MPC \\ \midrule
Jan	&	0.023	&	0.013	&	0.220	&	0.078	&	1.166	&	124.9	 &&	43.0	&	12.9	&	272.2	&	51.5	&	2.8	&		100\\
Feb	&	0.025	&	0.008	&	1.387	&	0.239	&	2.547	&	182.7	 &&	39.7	&	4.2	&	43.3	&	23.6	&	5.8	&		66\\
Mar	&	0.022	&	0.020	&	0.188	&	0.079	&	1.393	&	97.8	 &&	44.7	&	12.7	&	319.8	&	43.1	&	2.4	&		66\\
Apr	&	0.024	&	0.033	&	0.209	&	0.094	&	1.684	&	77.3	 &&	41.0	&	10.7	&	287.6	&	35.6	&	2.2	&		43\\
May	&	0.028	&	0.038	&	1.452	&	0.420	&	5.696	&	144.1	 &&	35.8	&	2.4	&	41.3	&	10.5	&	3.5	&	24\\
Jun	&	0.033	&	0.086	&	1.237	&	0.488	&	6.101	&	100.8	 &&	30.0	&	2.1	&	48.5	&	9.8	&	2.5	&		15\\
Jul	&	0.032	&	0.141	&	0.707	&	0.423	&	6.558	&	100.9	 &&	30.8	&	2.4	&	84.9	&	9.1	&	1.7	&		14\\
Aug	&	0.031	&	0.095	&	1.042	&	0.477	&	6.023	&	103.2	 &&	32.2	&	2.1	&	57.6	&	10.0	&	2.2	&		16\\
Sep	&	0.027	&	0.068	&	1.355	&	0.442	&	5.826	&	105.3	 &&	37.1	&	2.3	&	44.3	&	10.3	&	3.1	&		17\\
Oct	&	0.021	&	0.022	&	1.652	&	0.282	&	4.758	&	145.8	 &&	46.5	&	3.6	&	36.3	&	12.6	&	5.9			&	29\\
Nov	&	0.029	&	0.018	&	0.237	&	0.107	&	1.208	&	107.4	 &&	34.8	&	9.4	&	253.0	&	49.7	&	2.2	&		82\\
Dec	&	0.028	&	0.014	&	0.213	&	0.093	&	1.183	&	129.4	 &&	35.2	&	10.8	&	281.4	&	50.7	&	2.3	&		101\\
\bottomrule
\end{tabular}
\caption{Parameters for BLIP model, with independent within-cell pulse depths.}\label{par.BLIP}
\end{footnotesize}
\end{table}
\vspace{1cm}
\begin{table}[ht]
\begin{footnotesize}
\begin{tabular}{lrrrrrrrrrrrrrrr} \toprule
	&	$\lambda$	&	$\mu_X$	&	$\alpha$	&	$\alpha/\nu$	&	 $\kappa$	&	$\phi$	&	$\omega$	&&	MSIT	&	MSD	&	MCIT	&	 MCD	&	MCS	&	MPC\\ \midrule
Jan	&	0.024	&	0.001	&	2.147	&	4.591	&	1.027	&	0.046	&	 173	&&	42.3	&	8.9	&	23.8	&	24.5	&	22.4	&		165\\
Feb	&	0.023	&	0.001	&	3.680	&	4.394	&	1.096	&	0.058	&	 187	&&	42.6	&	5.4	&	17.1	&	18.8	&	18.8	&		177\\
Mar	&	0.023	&	0.001	&	2.000	&	5.525	&	0.712	&	0.043	&	 204	&&	44.1	&	8.3	&	30.5	&	21.7	&	16.4	&	 	195\\
Apr	&	0.024	&	0.001	&	2.000	&	6.740	&	0.517	&	0.039	&	 248	&&	41.7	&	7.7	&	34.4	&	17.8	&	13.4	&	 	239\\
May	&	0.027	&	0.001	&	2.000	&	7.760	&	0.437	&	0.054	&	 413	&&	37.3	&	4.8	&	35.4	&	15.5	&	8.1	&		392\\
Jun	&	0.031	&	0.001	&	2.000	&	9.607	&	0.310	&	0.050	&	 606	&&	32.1	&	4.1	&	40.3	&	12.5	&	6.2	&		576\\
Jul	&	0.030	&	0.001	&	2.000	&	10.413	&	0.167	&	0.039	&	 908	&&	33.4	&	4.9	&	69.2	&	11.5	&	4.2	&		874\\
Aug	&	0.029	&	0.001	&	2.000	&	9.683	&	0.293	&	0.053	&	 663	&&	34.4	&	3.9	&	42.2	&	12.4	&	5.6	&	 	630\\
Sep	&	0.025	&	0.001	&	2.000	&	8.901	&	0.345	&	0.047	&	 534	&&	40.1	&	4.8	&	39.1	&	13.5	&	7.4	&		510\\
Oct	&	0.021	&	0.001	&	2.126	&	6.698	&	0.580	&	0.041	&	 286	&&	48.4	&	6.9	&	29.2	&	16.9	&	14.3	&		274\\
Nov	&	0.025	&	0.001	&	2.000	&	5.389	&	1.055	&	0.049	&	 182	&&	39.9	&	7.6	&	21.1	&	22.3	&	21.5	&		173\\
Dec	&	0.026	&	0.001	&	2.035	&	4.584	&	1.093	&	0.054	&	 188	&&	37.9	&	7.9	&	23.6	&	25.7	&	20.1	&	 	179\\
\bottomrule
\end{tabular}
\caption{Parameters for BLIPR model, with common within-cell pulse depths; constraints: $\alpha > 2$, $\mu_X = 001$.}\label{par.BLIPRd2}
\end{footnotesize}
\end{table}
\vspace{1cm}
\begin{table}[ht]
\begin{footnotesize}
\begin{tabular}{lrrrrrrrrrrrrr} \toprule
	&	$\lambda$	&	$\iota$	&	$\alpha$	&	$\alpha/\nu$	&	 $\kappa$	&	$\phi$	&&	MSIT	&	MSD	&	MCIT	&	MCD	&	MCS		 \\ \midrule
Jan	&	0.022	&	0.164	&	2.075	&	5.014	&	0.996	&	0.042	 &&	46.2	&	9.1	&	23.2	&	23.1	&	24.6\\
Feb	&	0.021	&	0.177	&	3.451	&	4.818	&	1.063	&	0.053	 &&	47.5	&	5.5	&	16.5	&	17.5	&	20.9\\
Mar	&	0.020	&	0.196	&	2.000	&	5.910	&	0.695	&	0.041	 &&	48.8	&	8.3	&	29.2	&	20.3	&	18.0\\
Apr	&	0.022	&	0.241	&	2.000	&	7.083	&	0.509	&	0.037	 &&	46.5	&	7.6	&	33.3	&	16.9	&	14.8\\
May	&	0.023	&	0.400	&	2.000	&	8.127	&	0.434	&	0.052	 &&	43.9	&	4.7	&	34.0	&	14.8	&	9.4\\
Jun	&	0.026	&	0.586	&	2.000	&	10.015	&	0.311	&	0.049	 &&	38.9	&	4.1	&	38.5	&	12.0	&	7.3\\
Jul	&	0.024	&	0.879	&	2.000	&	10.777	&	0.173	&	0.040	 &&	42.3	&	4.6	&	64.3	&	11.1	&	5.3\\
Aug	&	0.024	&	0.639	&	2.000	&	10.109	&	0.299	&	0.052	 &&	42.3	&	3.8	&	39.7	&	11.9	&	6.8\\
Sep	&	0.021	&	0.518	&	2.000	&	9.257	&	0.343	&	0.045	 &&	47.4	&	4.8	&	37.7	&	13.0	&	8.6\\
Oct	&	0.019	&	0.277	&	2.051	&	7.006	&	0.575	&	0.039	 &&	53.8	&	7.1	&	29.1	&	16.7	&	15.7\\
Nov	&	0.023	&	0.175	&	2.000	&	5.832	&	1.018	&	0.045	 &&	43.9	&	7.6	&	20.2	&	20.6	&	23.5\\
Dec	&	0.024	&	0.179	&	2.000	&	5.018	&	1.056	&	0.050	 &&	42.0	&	8.0	&	22.6	&	23.9	&	22.2\\
\bottomrule
\end{tabular}
\caption{Parameters for $\text{BLRPR}_X$ model, with dependent intensity/duration ($\mu_X \propto \eta$); constraint: $\alpha > 2$.}\label{par.BLRPRM2}
\end{footnotesize}
\end{table}
\clearpage

 \newgeometry{bottom=2cm}
\section{Figures}\label{plots}
\begin{figure}[htpb]
\centering
  \includegraphics[width = 12cm]{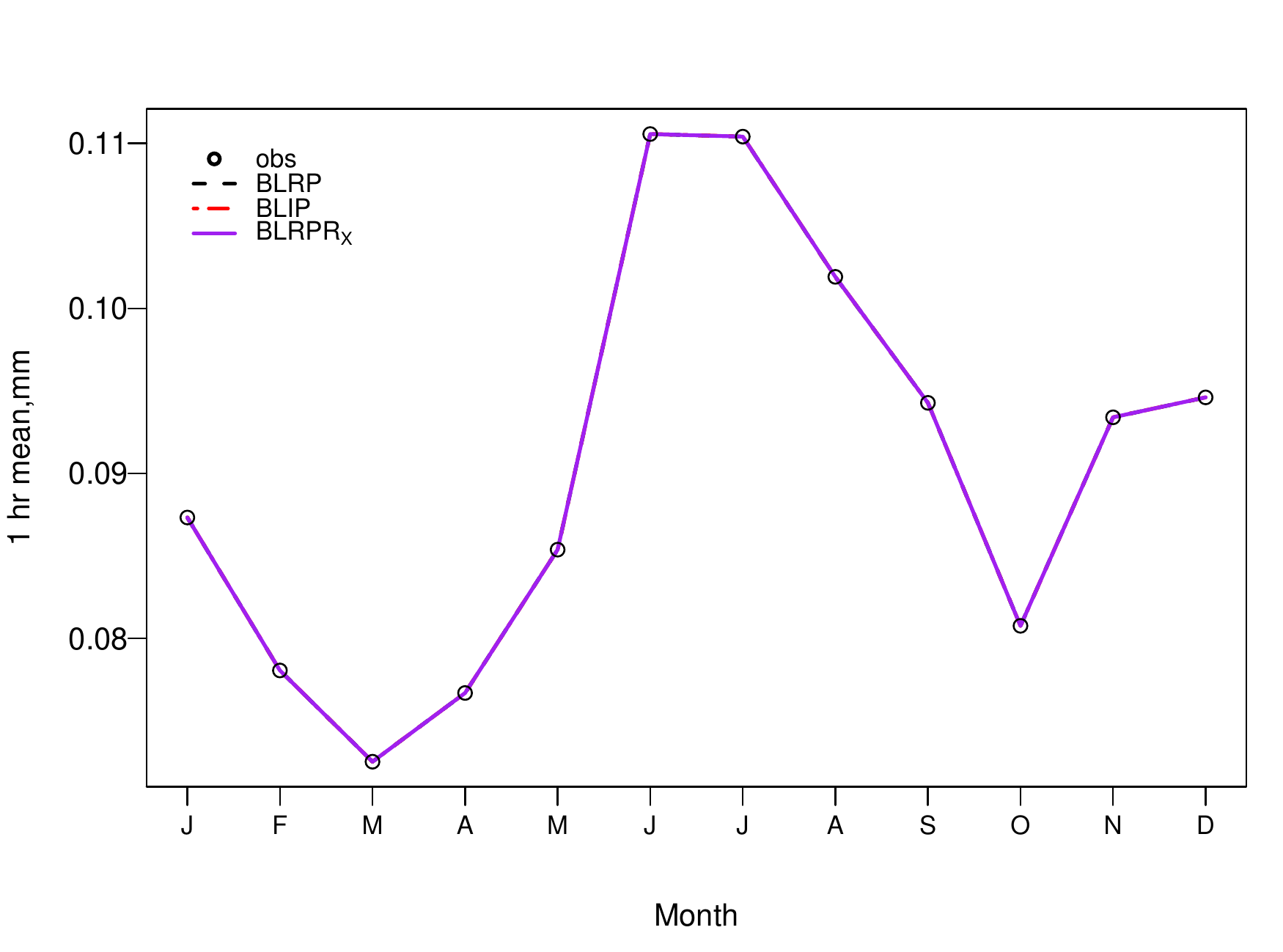}
  \caption{ Mean 1-hour rainfall by month, fitted v observed.}\label{mean}
   \vspace{20pt}
  \includegraphics[width = 13cm]{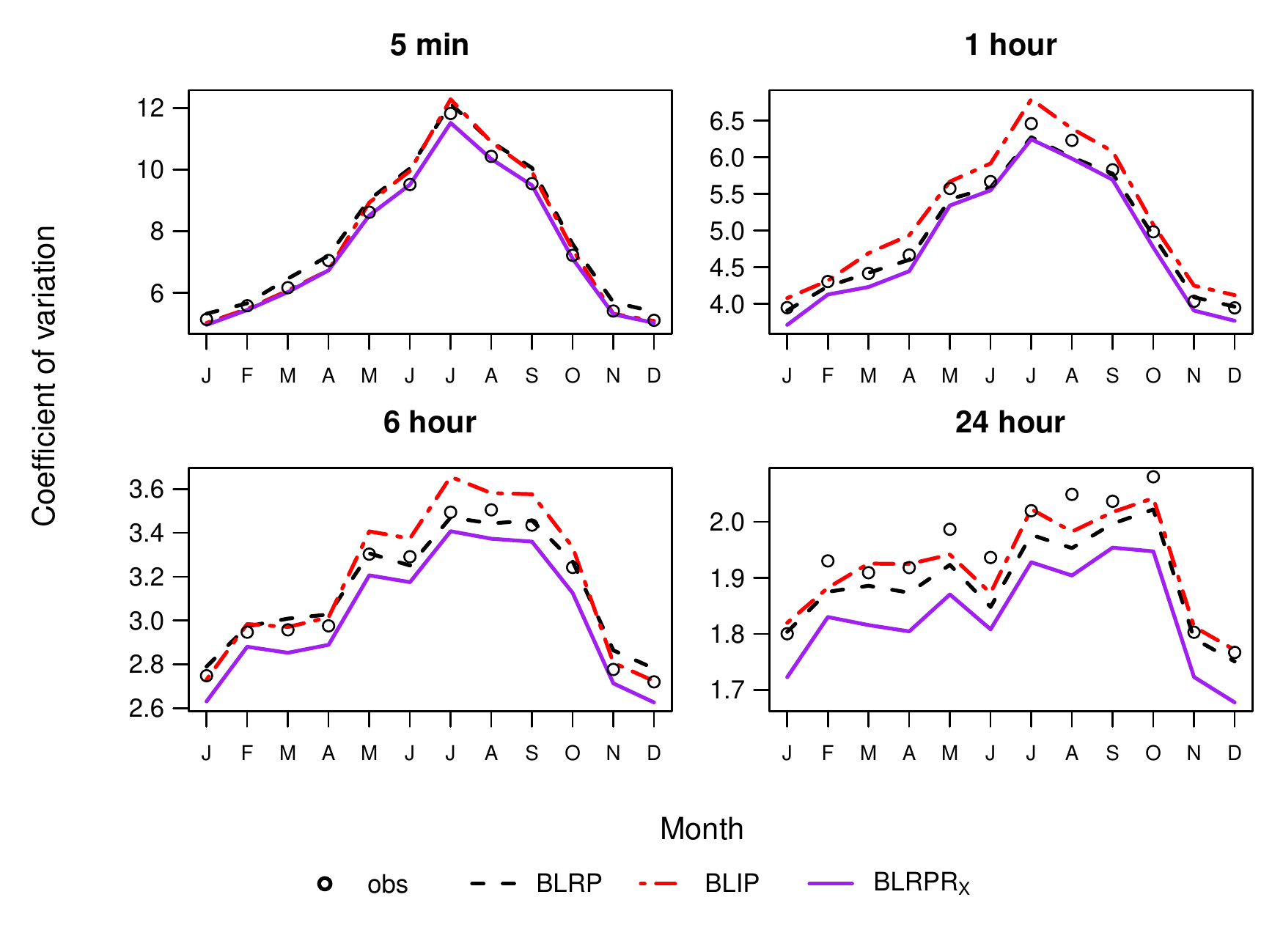}
  \caption{Coefficient of variation by month, fitted v observed.}\label{vars}
\end{figure}

\begin{figure}[htpb]
\centering
  \includegraphics[width = 13cm]{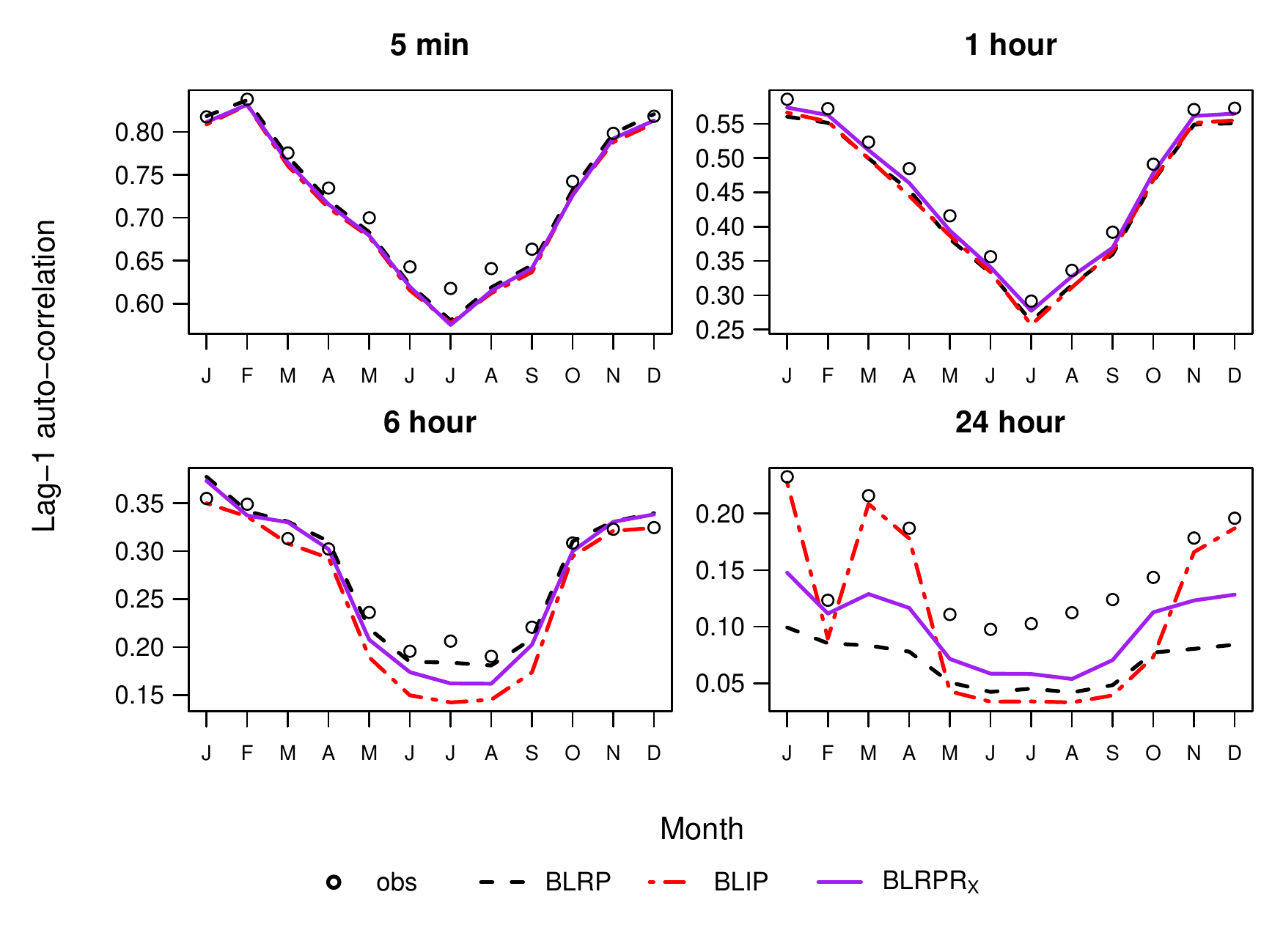}
 \caption{Lag-1 autocorrelation by month, fitted v observed.}\label{ac1}
  \includegraphics[width = 13cm]{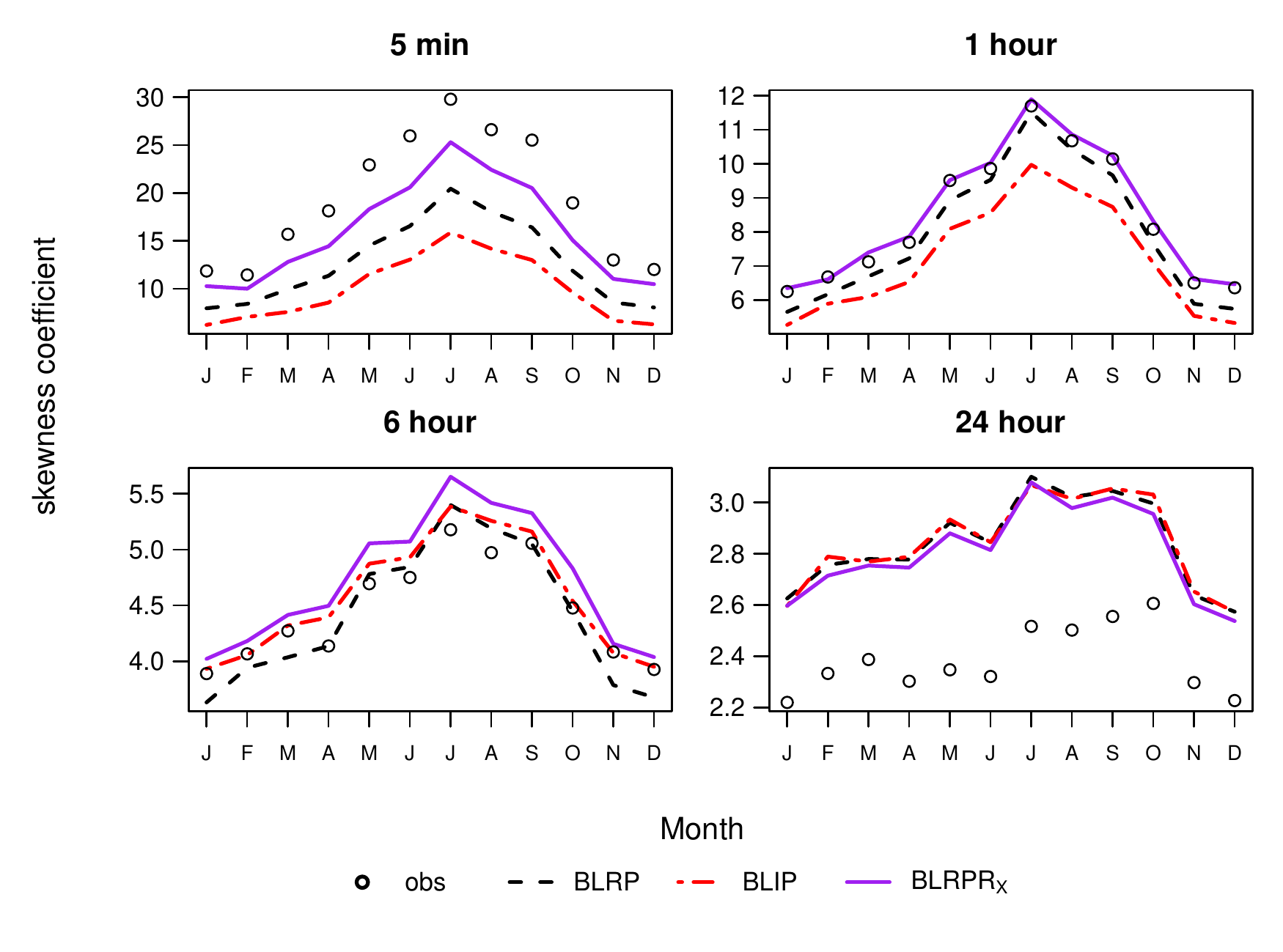}
  \caption{Coefficient of skewness by month, fitted v observed.}\label{skew}
\end{figure}
\begin{figure}[htpb]
\begin{center}
 \includegraphics[width = 13cm]{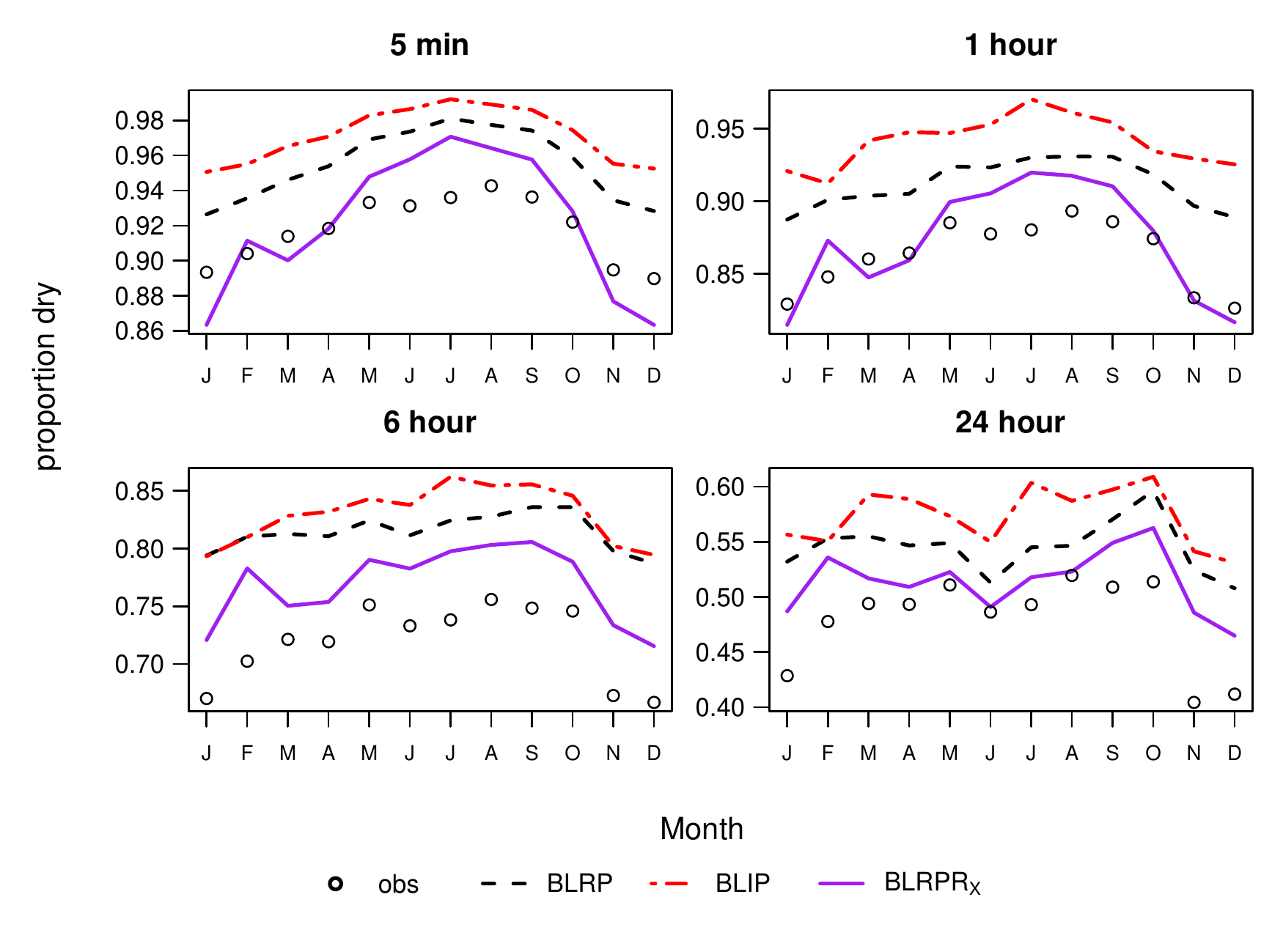}
  \caption{Proportion \add[JK]{of intervals that are} dry by month, fitted v observed.}\label{pdry}
    \vspace{20pt}
  \includegraphics[width = 13cm]{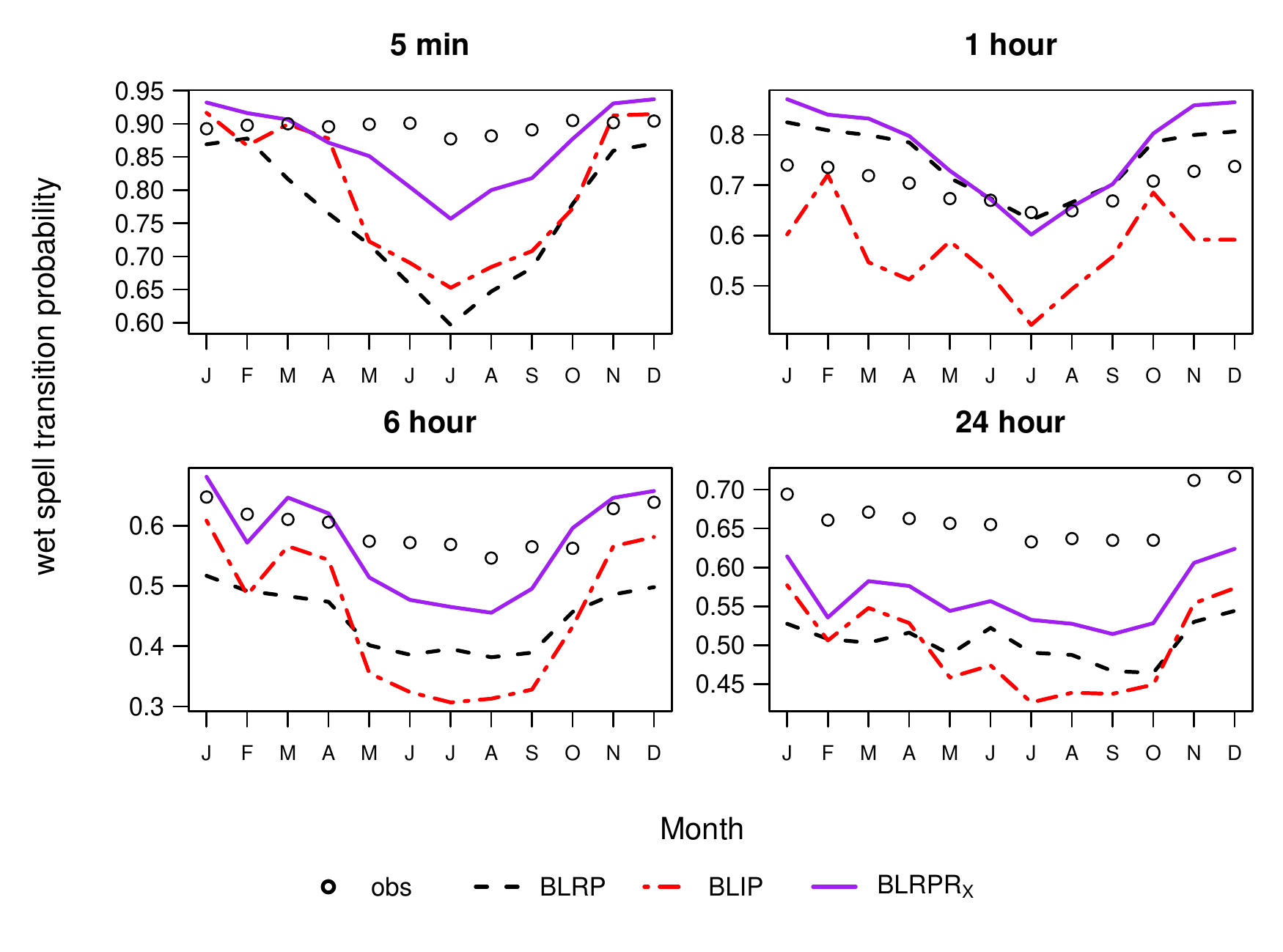}
  \caption{Transition probability of a wet interval being followed by another wet interval,  by month, fitted v observed.}\label{ww}
 \end{center}
\end{figure}
\begin{figure}[htp]
\centering
  \includegraphics[width = 12cm]{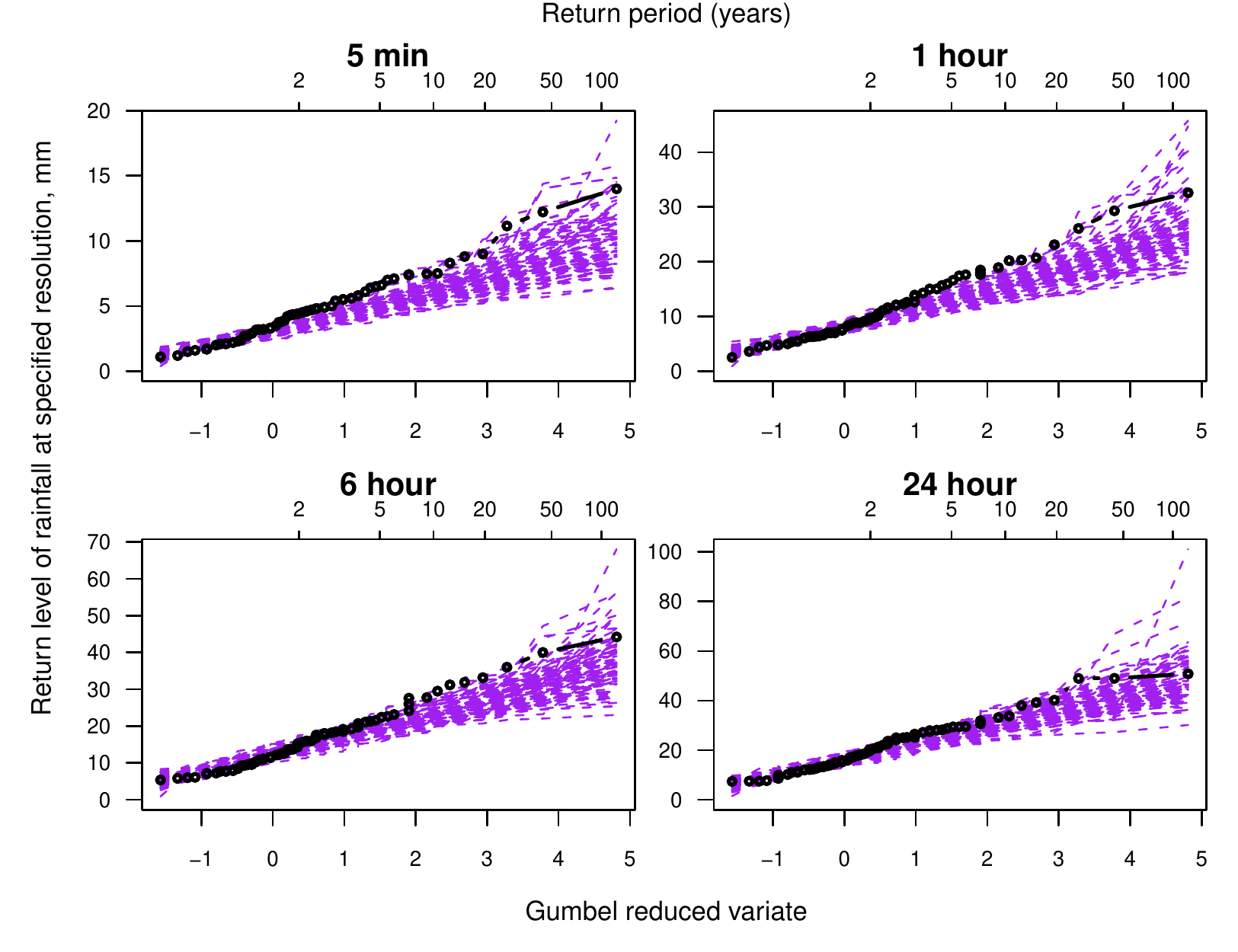}
\caption{Gumbel plots of observed (black) v simulated (purple) extremes for July, using the $\text{BLRPR}_X$ model and 100 simulations, each of 69 years; $\alpha$ constrained to be greater than 2.}\label{ext.jul}
\end{figure}
\begin{figure}[htpb]
\centering
\subfloat[5 minute]{\includegraphics[width = 7.9cm]{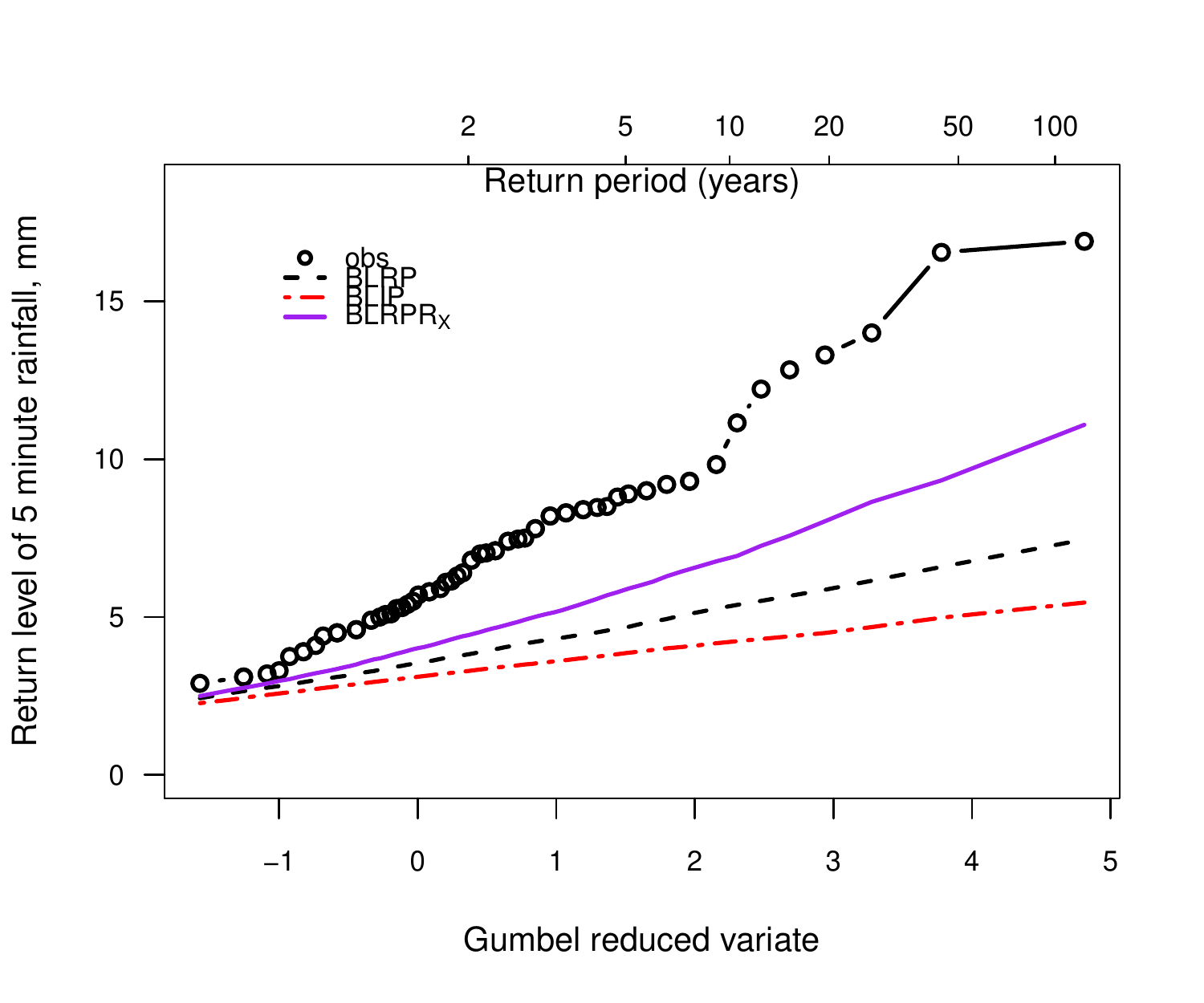}}
\subfloat[1 hour]{\includegraphics[width = 7.9cm]{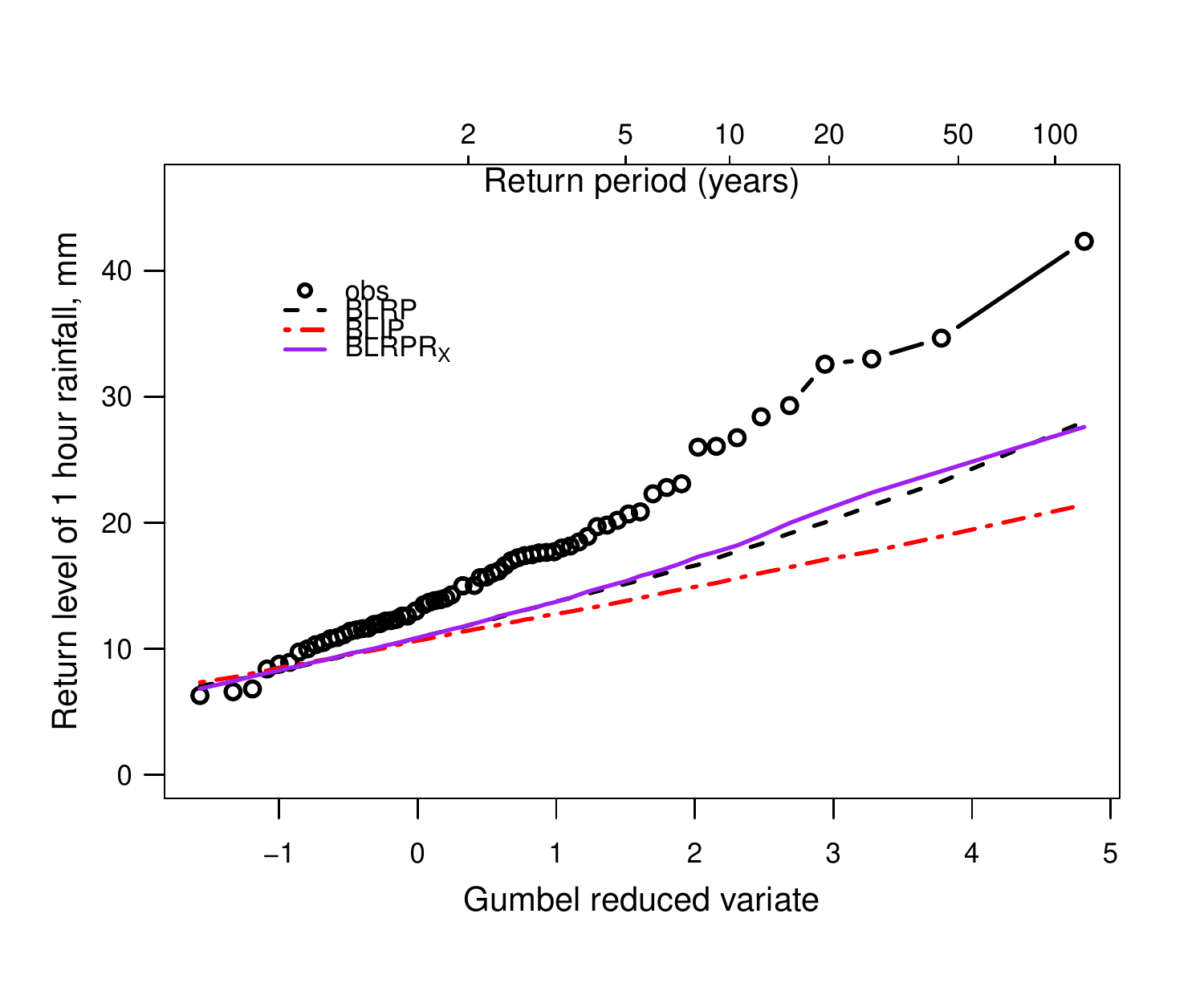}}
\caption{Annual Gumbel plots of observed v simulated extremes for variants of the Bartlett-Lewis model.}\label{ext.mean}
\end{figure}
\clearpage
\begin{figure}[htpb]
\centering
\subfloat[wide parameter range]{\includegraphics[width = 12.5cm]{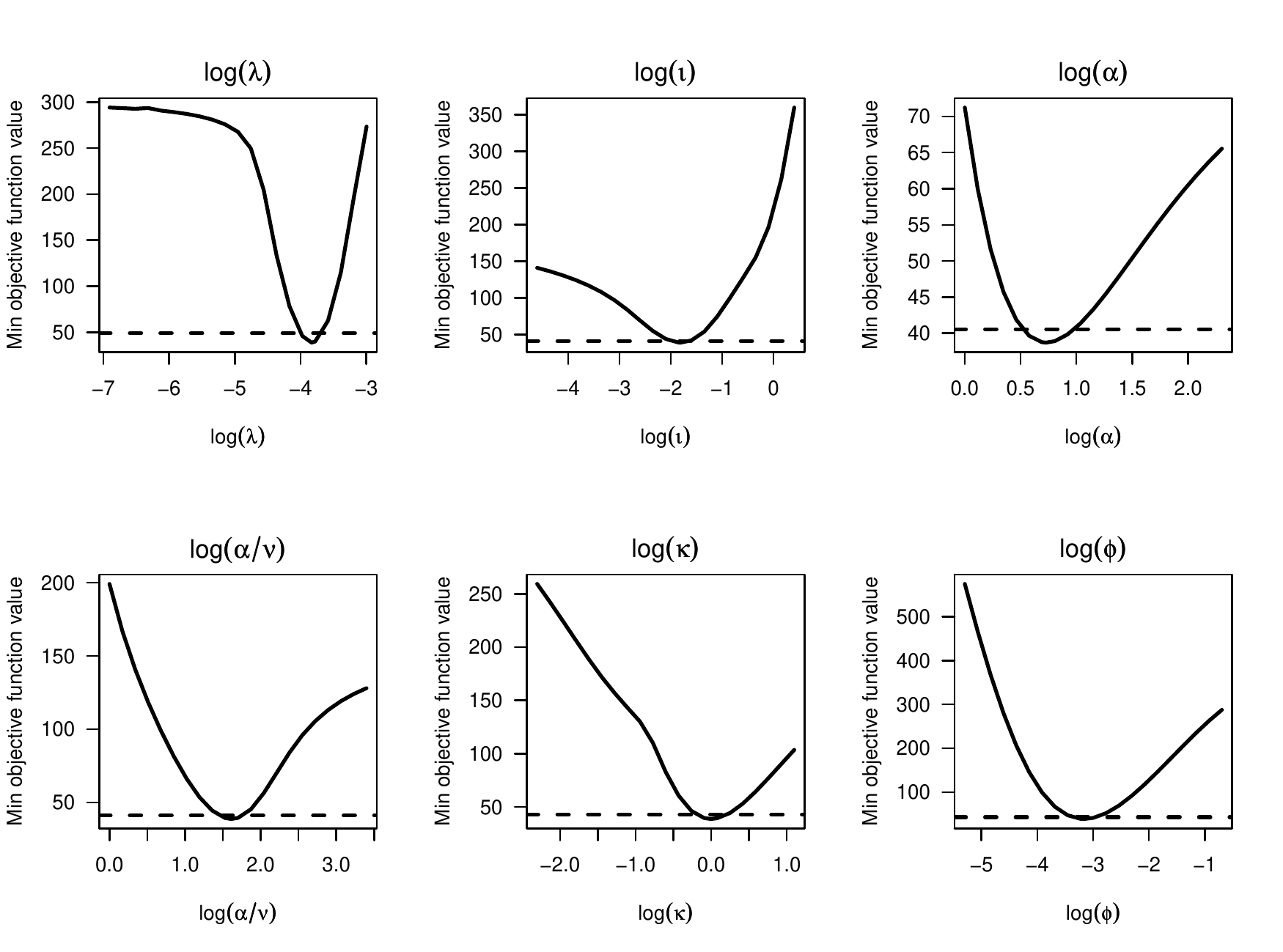}}\\
\subfloat[reduced parameter range]{\includegraphics[width = 12.5cm]{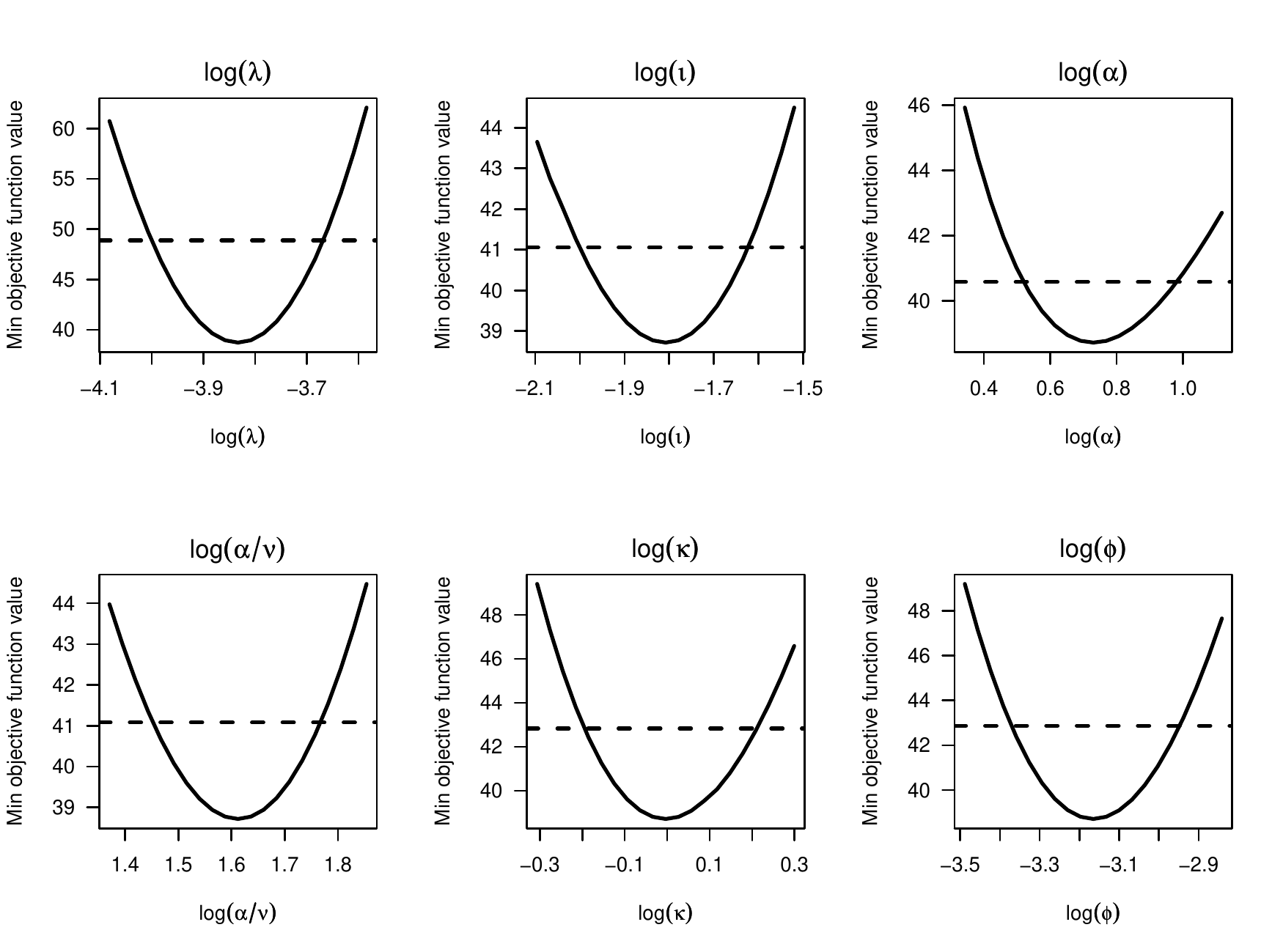}}
\caption{Profile objective function plots for the $\text{BLRPR}_X$ model for January; the plots show the logarithms of the parameters.}\label{profiles}
\clearpage
\end{figure}
\restoregeometry
\clearpage
\pagebreak
\end{appendix}
\bibliography{rainfall}
\bibliographystyle{agsm}
\end{document}